\documentclass[11pt,leqno]{article}

\usepackage{amssymb,amsmath,epsfig,amsthm}
\usepackage[mathscr]{eucal}
\input epsf
\newcommand{\n}{\noindent}

\newcommand{\intl}{\int\limits}
\newcommand{\vp}{\varepsilon}
\newcommand{\ovl}{\overline}
\newcommand{\bb}[1]{\mathbb{#1}}

\newcommand{\brel}[1]{\overset{\rightharpoonup}{#1}}

\numberwithin{equation}{section}

\theoremstyle{plain}
\newtheorem{thm}{Theorem}[section]
\newtheorem{lem}[thm]{Lemma}
\newtheorem{cor}[thm]{Corollary}

\theoremstyle{definition}

\newtheorem*{rk}{Remark}
\begin{document}
\begin{titlepage}
\title{\bf The exit distribution  for iterated Brownian motion in cones}
\author{Rodrigo Ba\~nuelos\thanks{Supported in part by NSF Grant
\# 9700585-DMS}\\Department of  Mathematics
\\Purdue University\\West Lafayette, IN
47906\\banuelos@math.purdue.edu\and Dante DeBlassie\\ Department of  Mathematics\\
Texas A\&M University\\College Station, TX \ 77843-3368\\
deblass@math.tamu.edu} \maketitle
\begin{abstract}
{\it We study the distribution of the exit place of iterated
Brownian motion in a cone, obtaining information about the chance
of the exit place having large magnitude. Along the way, we determine 
the joint distribution of the exit time and exit place of Brownian 
motion in a cone. This  yields information on large values of 
the exit place (harmonic measure) for Brownian motion. The harmonic measure 
for cones has been  studied by many authors for many years. Our results are 
sharper than any  previously obtained.}

\end{abstract}

\end{titlepage}

\section{Introduction}\label{sec1}

\indent

Roughly speaking, iterated Brownian motion (IBM) is ``Brownian motion run at an 
independent one-dimensional Brownian clock.'' Of course, this is not rigorous 
because the one-dimensional Brownian motion can take negative values, whereas 
Brownian motion is defined only for nonnegative times. There are two natural 
ways to get around this. First, one can use the absolute value of the 
one-dimensional Brownian motion. This process is one of the subjects of the 
papers Allouba and Zheng (2001) and Allouba (2002), where various connections 
with the biharmonic operator are presented. Those authors call their process 
``Brownian-time Brownian motion'' (BTBM). The other rigorous definition of IBM
is the one we will use and it 
is due to Burdzy (1993). He uses a natural extension of Brownian motion to 
negative times, called ``two-sided Brownian motion.'' 
Formally, let $X^+,X^-$ be independent 
$n$-dimensional Brownian motions started at $z\in {\bb R}^n$ and suppose $Y$ is 
one-dimensional Brownian motion started at 0, independent of $X^\pm$. Define 
two-sided Brownian motion by
\[
X(t) = \left\{\begin{array}{ll} X^+(t),&t\ge 0,\\
X^-(-t),&t<0.\end{array}\right.
\]
Then iterated Brownian motion is
\[
Z_t = X(Y_t),\quad t\ge 0.
\]

Although IBM is not a Markov process, it has many properties analogous to those 
of Brownian motion; we list a few here.  
\begin{enumerate}
\item[(1)] For instance, the process scales.  That is, for each $c>0$,
\[
cZ(c^{-4}t)
\]
is IBM.  
\item[(2)] The law of the iterated logarithm holds (Burdzy (1993))
\[
\limsup_{t\to 0} \frac{Z(t)}{t^{1/4} (\log \log (1/t))^{3/4}} = 
\frac{2^{5/4}}{3^{3/4}} \quad \text{a.s.}
\]
There is also a Chung--type LIL (Khoshnevisan and Lewis (1996)) and various Kesten--type
LIL's (Cs\"org\H{o}, F\"oldes and R\'ev\'esz (1996)) for IBM.  Other properties for local times
are proved in Xiao (1998). 
\item[(3)] The process has $4^{\rm th}$ order variation (Burdzy (1994)):
\[
\lim_{|\Lambda|\to 0} \sum^n_{k=1} [Z(t_k)-Z(t_{k-1})]^4 = 3(t-s) \text{ in } 
L^p,
\]
where $\Lambda = \{s=t_0 \le t_1 \le\cdots\le t_n = t\}$ is a partition of 
$[s,t]$ and $|\Lambda| = \max\limits_{1\le k\le n} |t_k-t_{k-1}|$. 
\end{enumerate}
An  interesting interpretation of IBM, due to Burdzy and Khoshnevisan (1998), is
as  a model for diffusion in a crack. See DeBlassie (2004) for other references.

There is a very interesting connection between IBM (as well as the BTBM process
of Allouba and Zheng) and the biharmonic operator $\Delta^2$. Namely, the function
\[u(t,x)=E_x[f(Z_t)]\]
solves the Cauchy problem
\begin{align*}
\frac{\partial}{\partial t}u(t,x)&=\frac{\Delta f(x)}{\sqrt{2\pi t}}+\frac 12\Delta^2u(t,x)\\
u(0,x)&=f(x)
\end{align*}
(Allouba and Zheng (2001) and DeBlassie (2004)). The appearance of the initial function
$f(x)$ in the PDE can be viewed as a manifestation of the non-Markovian nature of IBM.

This connection suggests the possibility of a relationship between IBM and 
initial-boundary or boundary value problems involving 
the biharmonic operator. While the results of DeBlassie (2004) are not encouraging for
connections with initial-boundary value problems, the work of Allouba and Zheng (2001)
suggests there is some hope for finding connections between probability and  Dirichlet-type boundary
value problems for the bilaplacian.  Such a connection, if found, would be particularly
exciting in its possible applications to the spectral theory (the study of eigenvalues and
eigenfunctions) of the bilaplacian where very little seems to be known.  An important first
step  in exploring this possibility, as in the case of the Laplacian and Brownian motion, is to
gain an understanding of the structure of the distribution of the exit place of IBM from open
sets, what  one may call, by abuse of terminology, the ``harmonic measure'' associated with
IBM.   In contrast with the BTBM process of Allouba and Zheng, this distribution does not
coincide with the usual harmonic measure associated with the Laplacian. The goal of this
article is to study the exit distribution of IBM from a cone in 
${\bb R}^n$. We chose this domain because it is unbounded and it contains a
boundary singularity. In addition, in this setting we are able to obtain explicit formulas
which 
lead to very sharp results. Our methods are easily adapted to bounded domains but in 
general our
formulas will not be as explicit and the result will not be as sharp.  

Let
$S^{n-1}$ be the unit sphere in
${\bb R}^n$. If $D$ is a proper  open subset of $S^{n-1}$, then the {\em generalized cone\/}
$C$ generated by $D$  is the set of rays emanating from the origin 0 passing through $D$.
Throughout  we assume $\partial D$ is $C^{2,\alpha}$. Then the Laplace--Beltrami operator 
$\Delta_{S^{n-1}}$ on $S^{n-1}$ with Dirichlet boundary conditions on $\partial 
D$ has a complete set of orthonormal eigenfunctions $m_j$ with corresponding 
eigenvalues $0<\lambda_1 < \lambda_2 \le \lambda_3\le\cdots$ such that
\begin{equation}\label{eq1.1}
\left\{\begin{array}{l}
\Delta_{S^{n-1}} m_j = -\lambda_j m_j\quad \text{and}\\
m_j = 0\quad \text{on}\quad \partial D,\end{array}\right.
\end{equation}
(Chavel (1984)). If $B$ is $n$-dimensional Brownian motion and $\tau_C(B)$ is 
its exit time from $C$, then it is known (DeBlassie (1988))
\begin{equation}\label{eq1.2}
P_x(\tau_C(B) > t) \sim C(x) t^{-p_1/2}\quad \text{as}\quad t\to \infty
\end{equation}
where
\begin{equation}\label{eq1.3}
p_1 = \sqrt{\lambda_1 + \left(\frac{n}{2}-1\right)^2} - \left(\frac{n}2 - 
1\right)
\end{equation}
and
\begin{equation}\label{eq1.4}
C(x) = \left(\frac{|x|^2}2\right)^{p_1/2} \frac{\Gamma\left(\frac12(p_1 + 
n)\right)}{\Gamma\left(p_1+\frac{n}2\right)} \left[~\intl_{\partial D} 
m_1(\eta) 
\mu(d\eta)\right] m_1\left(\frac{x}{|x|}\right),
\end{equation}
$\mu$ being surface measure on $S^{n-1}$. Here and in what follows,
\[
f(t) \sim g(t) \quad \text{as}\quad t\to\infty
\]
means
\[
\frac{f(t)}{g(t)}\to 1 \quad \text{as}\quad t\to\infty.
\]

In DeBlassie (2004) it is shown that if $\tau_C(Z)$ is the first exit time of 
IBM $Z$ from $C$, then as $t\to \infty$,
\begin{equation}\label{eq1.5}
P_x(\tau_C(Z) > t) \approx \begin{cases}
t^{-p_1/2},&\text{$p_1<2$}\\
t^{-1} \ln t,&\text{$p_1=2$}\\
t^{-(p_1+1)/2},&\text{$p_1>2$}\end{cases}
\end{equation}
where $f(t) \approx g(t)$ means there exist constants $C_1$ and $C_2$ such that 
\[
C_1 \le \frac{f(t)}{g(t)} \le C_2,\quad t \text{ large.}
\]
In light of Burkholder's inequalities (1977) and \eqref{eq1.2},
\begin{equation}\label{eq1.6}
E_x(|B(\tau_C)|^p) < \infty\quad \text{iff}\quad p<p_1.
\end{equation}
Hence considering the ``fourth order'' properties of IBM described above, we 
expect \eqref{eq1.5} should imply
\[
E_x(|Z(\tau_C)|^p] < \infty\quad \text{iff}\quad \begin{cases}
p<2p_1,&\text{$p_1\le 2$}\\
p<2(p_1+2),&\text{$p_1>2$.}\end{cases}
\]
Indeed, we have the following theorem. We will always assume the positive 
$x_n$-axis passes through $C$. If $\varphi(\eta)$ is the angle between $\eta\in 
S^{n-1}$ and the positive $x_n$-axis, then in polar coordinates $y=r\eta$, the 
$(n-1)$-dimensional surface measure $\sigma$ on $\partial C$ is given by
\begin{equation}\label{eq1.7}
\sigma(dy) = r^{n-2} \sin \varphi(\eta) \mu(d\eta) dr.
\end{equation}

\begin{thm}\label{thm1.1}
As $r\to \infty$, for $z=\rho\theta$,
\[
\frac{d}{dr} P_z(|Z(\tau_C)|\le r) \sim A(z,p_1) \begin{cases}
r^{-2p_1-1},&\text{$\dfrac{p_1}2 < 1$},\\
\noalign{\smallskip}
r^{-5}\ln r,&\text{$\dfrac{p_1}2 =1$},\\
\noalign{\smallskip}
r^{-p_1-3},&\text{$\dfrac{p_1}2 >1$,}
\end{cases}
\]
where for $\frac{p_1}2 < 1$,
\begin{eqnarray*}
A(z,p_1) = \rho^{2p_1} m^2_1(\theta) \frac{\Gamma\left(\frac{p_1+n}2\right) 
\Gamma\left(\frac{3p_1+n}2-1\right)}{\left(\Gamma\left(p_1 + 
\frac{n}2\right)\right)^2} \left[~\intl_{\partial D} m_1(\eta)\mu(d\eta)\right]\cdot\\ 
\cdot\left[~\intl_{\partial D} \sin \varphi(\eta)
\frac\partial{\partial n_\eta} m_1(\eta)\mu(d\eta)\right] \int^\infty_0 w^{-p_1/2} (1+w)^2
dw,
\end{eqnarray*}
the integrals over $\partial D$ are taken with respect to $\mu(d\eta)$ and 
$\frac\partial{\partial n_\eta}$ denotes the inward normal derivative at 
$\partial D$; for $\frac{p_1}2 = 1$,
\[
A(z,p_1) = 2\left(1 + \frac{n}2\right)^{-1} \rho^4 m^2_1(\theta) 
\left[~\intl_{\partial D} m_1(\eta)\mu(d\eta)\right] \left[~\intl_{\partial D}
 \sin \varphi(\eta)
\frac\partial{\partial n_\eta} m_1(\eta)\mu(d\eta)\right];
\]
and for $\frac{p_1}2 > 1$,
\[
A(z,p_1) = 2\rho^{p_1} m_1(\theta) \left[~\intl_{\partial D} \sin(\eta)\varphi 
\frac\partial{\partial n_\eta} m_1(\eta)\mu(d\eta)\right] E_z(\tau_{BM}),
\]
where $\tau_{BM}$ is the first exit time of Brownian motion from $C$.
\end{thm}

\begin{cor}\label{cor1.2}
a)~~As $r\to \infty$,
\[
P_z(|Z(\tau_C)| > r) \sim A(z,p_1) \begin{cases}
\dfrac{1}{2p_1} r^{-2p_1},&\text{$\dfrac{p_1}2 < 1$},\\
\noalign{\smallskip}
\dfrac14 r^{-4} \ln r,&\text{$\dfrac{p_1}2 =1$,}\\
\noalign{\smallskip}
\dfrac1{p_1+2} r^{-p_1-2},&\text{$\dfrac{p_1}2>1$.}
\end{cases}
\]
b)~~We have
\[
E_z[|Z(\tau_C)|^p] < \infty\quad \text{iff}\quad \begin{cases}
p<2p_1,&\text{$\dfrac{p_1}2<1$,}\\
\noalign{\smallskip}
p<4,&\text{$\dfrac{p_1}2 = 1$,}\\
\noalign{\smallskip}
p<p_1+2,&\text{$\dfrac{p_1}2>1$,}
\end{cases}
\]
\end{cor}

\begin{rk}
Below in \eqref{eq3.3} we get a series expansion of the density $\frac{d}{dr} 
P_z(|Z(\tau_C)| \le r)$ valid for $r\ne |z|$, but it is not all that 
enlightening.
\end{rk}

\n Along the way to proving Theorem \ref{thm1.1}, we derive the following result 
of independent interest.

\begin{thm}\label{thm1.3}
For $\tau$ being the first exit time of Brownian motion $B$ from the cone $C$,
\[
P_x(B_\tau \in dy, \tau\in dt) = \frac12 \frac\partial{\partial n_y} 
p_C(t,x,y) \sigma(dy)dt,
\]
where $\frac\partial{\partial n_y}$ is the inward normal derivative at $\partial 
C$, $p_C(t,x,y)$ is the transition density of Brownian motion killed upon 
exiting $C$ and $\sigma$ is surface measure on $\partial C$.\hfill $\square$
\end{thm}

Hsu (1986) has proved this result for bounded $C^3$ domains. But because the 
cone $C$ is unbounded with a boundary singularity, there are technicalities not 
present in the case considered by Hsu.

We have the following consequence of Theorem \ref{thm1.3} that is also of 
independent interest. Note it gives an improvement of \eqref{eq1.6} above.

\begin{thm}\label{thm1.4}
Let $\tau$ be the exit time of Brownian motion from $C$. Then for $x=\rho\theta$ 
and $r\ne \rho$, 
\begin{eqnarray*}
\frac{d}{dr} P_x(|B_\tau|\le r) &= \frac12 r^{\frac{n}2-2} \rho^{1-\frac{n}2} 
\sum^\infty_{j=1} \alpha^{-1}_j \gamma^{\alpha_j} 
[1+(1-\gamma^2)^{1/2}]^{-\alpha_j} \cdot\\
&\quad\cdot\intl_{\partial D} \sin \varphi(\eta) 
\left[\frac\partial{\partial n_\eta} m_j(\eta)\right] \mu(d\eta) m_j(\theta),
\end{eqnarray*}
where $\alpha_j = \sqrt{\lambda_j - \left(\frac{n}2-1\right)^2}$ and $\gamma = 
\frac{2\rho r}{\rho^2 +r^2}$ and the convergence is uniform for $\gamma\le 
1-\vp$.\hfill $\square$
\end{thm}

\begin{cor}\label{cor1.5}
As $r\to \infty$, for $x=\rho\theta$,
\[
P_x(|B_\tau|>r) \sim \frac{2^{p_1+\frac{n}2-2} \rho^{p_1}}{p_1(p_1 + 
\frac{n}2-1)} 
\left[~\intl_{\partial D} \sin \varphi(\eta) \frac\partial{\partial n_\eta} 
m_1(\eta) \mu(d\eta)\right] m_1(\theta) r^{-p_1}.
\]
\end{cor}

It follows from the classical estimates for harmonic measure (see Haliste (1984) and 
Essen and Haliste (1984)) that there are constants $C_1$ and $C_2$, depending on $x$,  such
that for large
$r$, 
\[
C_1 r^{-p_1} \leq P_x(|B_\tau|>r)\leq C_2 r^{-p_1}.
\]
However, as far as we know these techniques do not identify the exact limit as Corollary
\ref{cor1.5} above does.  It is also interesting to note here that in the case of the
parabolic--shaped regions 
\begin{equation}
{\cal P}_\alpha=\{(x,Y)\in\bb R\times\bb R^{n-1}\colon x > 0,\, |Y| < Ax^\alpha\}\nonumber,
\label{parabola1}
\end{equation}
with  $0 < \alpha < 1$ and $A>0$, it is proved in Ba\~nuelos and Carroll (2003) that 
\begin{equation}
\log 
       P_{z}(|B_\tau |>r) \sim  -\frac{\sqrt{\lambda_1}}{A(1-\alpha)}
r^{1-\alpha}, 
\label{parabola2}
\end{equation}
where  $\lambda_1$ is the smallest eigenvalue for the Dirichlet
Laplacian in the unit ball of $\bb R^{n-1}$.  In view of Corollary \ref{cor1.5}, it is natural to ask
if it is possible to obtain a similar expression for the harmonic measure of the parabolic-shaped
regions, and in particular to identify the asymptotics of 
$P_{z}(|B_{\tau}|>r)$. That is, is it possible to obtain
a result similar to that in Ba\~nuelos and
Carroll but without the logs?   At present we do not know the answer to this question.  For
various results related to the asymptotics of  exit times of Bronian motion and heat kernels for
parabolic--type regions, we refer the reader to \cite{ba}, \cite{van}, \cite{DSmits}, \cite{li},
\cite{lif}. 
 
Finally, Allouba and Zheng (2001) show the exit distribution of their BTBM process is the 
same as that of Brownian motion---i.e., harmonic measure  (see their Theorem 0.2). In 
light of this, Theorem \ref{thm1.4} above yields the density of the size of the 
exit place of BTBM in a cone. Also note for IBM, the exit distribution is NOT 
the same as the exit distribution of two-sided Brownian motion in $C$.

The article is organized as follows. In section 2 we establish various estimates on
the terms in the series expansion of the heat kernel of the cone. Then we use them
to prove Theorem 1.4 and Corollary 1.5. Using Theorem 1.3, we prove Theorem 1.1 in
section 3. The proof of Theorem 1.3 is given in section 4, using some results of
Pinsky. The proof is independent of the proof of Theorem 1.1.

\section{Auxiliary Results; Proof of Theorem \ref{thm1.4} and Corollary 
\ref{cor1.5}}\label{sec2}

\indent

In what follows, we will make repeated use of the following result (Lemma 6.18 
on page 111) from Gilbarg and Trudinger (1983).

\medskip
\n {\bf Elliptic Regularity Theorem.} {\em Suppose
\[
L = \sum_{i,j} a_{ij}(x) \frac{\partial^2}{\partial x_i \partial x_j} + \sum_j 
b_j(x) \frac\partial{\partial x_j} + c(x)
\]
is a strictly elliptic operator on a domain $\Omega\subseteq {\bb R}^n$. Assume 
the coefficients of $L$ are in $C^\alpha(\ovl\Omega)$, $\Omega$ has a 
$C^{2,\alpha}$ boundary portion $T$ and $\varphi\in C^{2,\alpha}(\ovl\Omega)$. 
If $u\in C(\ovl\Omega) \cap C^2(\Omega)$ satisfies $Lu = 0$ in $\Omega$ and 
$u=\varphi$ on $T$, then $u\in C^{2,\alpha}(\Omega\cup T)$.$\hfill \square$
}\medskip

The heat kernel for $C$ has a series expansion, due to Ba\~nuelos and Smits 
(1997): \ For $x=\rho\theta$, $y=r\eta$, 
\begin{equation}\label{eq2.1}
p_C(t,x,y) = t^{-1} (\rho r)^{1-\frac{n}2} e^{-\frac{\rho^2+r^2}{2t}} 
\sum^\infty_{j=1} I_{\alpha_j}\left(\frac{\rho r}t\right) m_j(\theta) m_j(\eta),
\end{equation}
where $\rho,r>0$, $\theta,\eta \in S^{n-1}$, and
\begin{equation}\label{eq2.2}
\alpha_j = \sqrt{\lambda_j + \left(\frac{n}2-1\right)^2}.
\end{equation}
The convergence is uniform for $(t,x,y) \in (T,\infty)\times \{x\in C\colon \ 
|x|<R\} \times C$, for any positive constants $T$ and $R$.
 The 
modified Bessel function $I_\nu$ is given by
\[
I_\nu(z) = \left(\frac{z}2\right)^\nu \sum^\infty_{k=0} 
\left(\frac{z}2\right)^{2k} \frac1{k!\Gamma(\nu+k+1)}.
\]
First we show termwise normal differentiation at the boundary is permitted. Recall
$\frac{\partial}{\partial n_y}$ and $\frac{\partial}{\partial n_\eta}$ denote inward
normal differentiation at $y\in\partial C\backslash\{0\}$ and $\eta\in\partial D$, 
respectively.

\begin{lem}\label{lem2.1}
For  $y=r\eta\in \partial C\backslash\{0\}$,
\[
\frac\partial{\partial n_y} p_C(t,x,y) = r^{-1}t^{-1} (\rho r)^{1-\frac{n}2} 
e^{-\frac{\rho^2+r^2}{2t}} \sum^\infty_{j=1} I_{\alpha_j} \left(\frac{\rho 
r}t\right) m_j(\theta) \frac\partial{\partial n_\eta} m_j(\eta)
\]
uniformly for $(t,x,y) \in (T,\infty) \times \{x\in C\colon \ |x|<R\}$ 
$\times \{x\in \partial C\colon \ |x|<R\}$, where $T,R>0$ are arbitrary.
\end{lem}

\begin{proof}
We have
\[
\nabla = \brel e_r \frac\partial{\partial r} + \frac1r \nabla_{S^{n-1}},
\]
where $\brel e_r$ is a unit vector in the radial direction and 
$\nabla_{S^{n-1}}$ is the gradient operator on $S^{n-1}$. Thus for $y\in 
\partial C\backslash \{0\}$
\[
\frac\partial{\partial n_y} m_j(\eta) = \frac1r \frac\partial{\partial n_\eta} 
m_j(\eta).
\]
Consequently we need only verify the uniform convergence of
\begin{equation}\label{eq2.3}
\sum^\infty_{j=1} I_{\alpha_j}\left(\frac{\rho r}t\right) m_j(\theta) 
\frac\partial{\partial n_\eta} m_j(\eta)
\end{equation}
for $(t,x,y) \in B = (T,\infty) \times \{x\in C\colon \ |x|<R\}\times \{x\in 
\partial C\colon \ |x|<R\}$.

Since $\|m_j\|_2 = 1$, by Theorem 8 on page 102 of Chavel (1984), for some positive
$c_n$ and $b_n=b(n)$ depending only on $n$,
\begin{equation}\label{eq2.4}
\sup_D |m_j| \le c_n\lambda^{b(n)/4}_j.
\end{equation}

Since $(\Delta_{S^{n-1}} - \lambda_j)m_j = 0$ on $D$ and $m_j\in C(\ovl D)$, by 
the Elliptic Regularity Theorem, $m_j\in C^{2,\alpha}(\ovl D)$. Hence by the 
$C^{2,\alpha}$ nature of $\partial D$ and the global Schauder estimates (Theorem 
6.6 on page 98 of Gilbarg and Trudinger (1983)), for some constant $K$ 
independent of $j$,
\begin{align}
\sup_{\ovl D} |\nabla_{S^{n-1}}m_j| &\le K\sup_D |m_j|\nonumber\\
\label{eq2.5}
&\le K \lambda^{b(n)/4}_j.
\end{align}
Here and in what follows, $K$ will be a number whose value might change from 
line to line, but is independent of $j$. Hence for $\theta\in D$ and $\eta\in 
\partial D$ (using \eqref{eq2.2})
\begin{equation}\label{eq2.6}
\left|m_j(\theta) \frac\partial{\partial n_\eta} m_j(\eta)\right| \le 
K\alpha^{b(n)}_j.
\end{equation}
By formula \eqref{eq2.4} on page 303 in Ba\~nuelos and Smits (1997),
\begin{equation}\label{eq2.7}
I_\nu(z) \le K\left(\frac{z}2\right)^\nu \frac{e^\nu e^z}{\left(\nu + 
\frac12\right)^\nu}
\end{equation}
where $K$ is independent of $\nu$ and $z$. Then to show uniform convergence of 
\eqref{eq2.3} on $B$, it suffices to show for $M = \frac{R^2}T$,
\begin{equation}\label{eq2.8}
\sum^\infty_{j=1} \frac{(Me)^{\alpha_j} \alpha^{b(n)}_j}{\left(\alpha_j + 
\frac12\right)^{\alpha_j}} < \infty.
\end{equation}
By the Weyl asymptotic formula (Chavel (1984) page 172), there are constants 
$K_1$ and $K_2$ such that
\begin{equation}\label{eq2.9}
K_1j \le \alpha^{n-1}_j \le K_2j,\qquad j\ge 1.
\end{equation}
Then for $j$ large, 
\begin{align*}
\left(\alpha_j + \frac12\right)^{\alpha_j}&\ge (\alpha_j)^{\alpha_j}\\
&\ge \exp (K_1 j^{\frac1{n-1}} \ln (K_1 j^{\frac1{n-1}})).
\end{align*}
Hence for $c = K_1$ if $Me<1$ and $c = K_2$ if $Me\ge 1$, the sum in 
\eqref{eq2.8} is bounded by
\begin{align*}
&\sum^\infty_{j=1} j^{b(n)/(n-1)} \exp(cj^{\frac1{n-1}} \ln(Me) - K_1 j^{\frac1{n-1}} \ln 
(K_1j^{\frac1{n-1}}))\\
&\quad \le K \sum^\infty_{j=1} j^{b(n)/(n-1)} \exp \left(-\frac12 K_1 
j^{\frac1{n-1}}\right)\\
&\quad < \infty, \text{ by the integral test.}\qquad \qed
\end{align*}
\renewcommand{\qed}{}\end{proof}

\n In the sequel, we will use the following bound, which is an immediate 
consequence of \eqref{eq2.4} and \eqref{eq2.5}.

\begin{cor}\label{cor2.2}
For some positive $K$ and $b(n)$ independent of $j$,
$$
\sup_D |m_j| \vee \sup_{\ovl D} |\nabla_{S^{n-1}} m_j| \le K\alpha^{b(n)/2}_j. 
\eqno \square
$$
\end{cor}

\begin{lem}\label{lem2.3}
For some $K>0$, with $\gamma = \frac{2\rho r}{\rho^2+r^2}$,

a)~~$\int^\infty_0 t^{-1} e^{-\frac{\rho^2+r^2}{2t}} I_\alpha\left(\frac{\rho 
r}{t}\right)dt \le K\alpha^{-1} \gamma^\alpha$

b)~~$\left|\frac\partial{\partial \rho} \int^\infty_0 t^{-1} 
e^{-\frac{\rho^2+r^2}{2t}} I_\alpha \left(\frac{\rho r}t\right) dt\right| \le 
\frac{r(1-\gamma^2)^{-\frac12} (\gamma^{\alpha-1} + 
\gamma^{\alpha+1})}{\rho^2+r^2}$

\n for positive $\rho,r$ and $\alpha$, with $\rho\ne r$.
\end{lem}

\begin{proof}
Change variables $w = \frac{\rho^2+r^2}{2t}$ to get
\begin{align}
\int^\infty_0 t^{-1} e^{-\frac{\rho^2+r^2}{2t}} I_\alpha \left(\frac{\rho 
r}{t}\right) dt &= \int^\infty_0 \frac{2w}{\rho^2+r^2} e^{-w} I_\alpha 
\left(\frac{2\rho rw}{\rho^2+r^2}\right) \frac{\rho^2+r^2}{2w^2} dw\nonumber\\
\label{eq2.10}
&= \int^\infty_0 w^{-1}e^{-w}I_\alpha(\gamma w) dw.
\end{align}
Notice since $\rho \ne r, 0<\gamma<1$.

\n a)~~Using the expansion of $I_\nu(z)$ given before Lemma \ref{lem2.1}, by 
monotone convergence
\begin{align*}
\int^\infty_0 w^{-1}e^{-w} I_\alpha(\gamma w)dw &= \sum^\infty_{k=0} 
\left(\frac\gamma2\right)^{\alpha+2k} \frac1{k!\Gamma(\alpha+1+k)} \int^\infty_0 
w^{\alpha+2k-1} e^{-w}\ dw\\
&= \sum^\infty_{k=0} \left(\frac\gamma2\right)^{\alpha+2k} 
\frac{\Gamma(\alpha+2k)}{k!\Gamma(\alpha+1+k)}\\
&= \frac1{2\sqrt\pi} \sum^\infty_{k=0} \gamma^{\alpha+2k} 
\frac{\Gamma\left(\frac\alpha2 + k\right) \Gamma\left(\frac\alpha2 + 
k+\frac12\right)}{k!\Gamma(\alpha+1+k)},
\end{align*}
where we have used the formula
\begin{equation}\label{eq2.11}
\Gamma(2z) = \frac1{\sqrt{2\pi}} 2^{2z-\frac12} \Gamma(z) \Gamma\left(z + 
\frac12\right)
\end{equation}
(Abramowitz and Stegun (1972) page 256, 6.1.18) for $z = \frac\alpha2 +k$. We 
also make use of the following formulas from Abramowitz and Stegun (1972) for 
the hypergeometric function $F$:
\begin{align*}
&F(a,b;c,z) = \frac{\Gamma(c)}{\Gamma(a)\Gamma(b)} \sum^\infty_{n=0} 
\frac{\Gamma(a+n) \Gamma(b+n)}{\Gamma(c+n)} \frac{z^n}{n!}\\
&F\left(a,\frac12+a; 1+2a,z\right) = 2^{2a} [1+(1-z)^{1/2}]^{-2a}
\end{align*}
(formulas 15.1.1 and 15.1.13 on page 556). Using the first one, then the second, 
yields
\begin{align}
\int^\infty_0 w^{-1}e^{-w} I_\alpha(\gamma w)dw &= 
\frac{\gamma^\alpha}{2\sqrt\pi} \frac{\Gamma\left(\frac\alpha2\right) 
\Gamma\left(\frac{\alpha+1}2\right)}{\Gamma(\alpha+1)} 
F\left(\frac\alpha2, \frac\alpha2 + \frac12; \alpha+1, 
\gamma^2\right)\nonumber\\
&= \frac{\gamma^\alpha}{2\sqrt\pi} \frac{\Gamma\left(\frac\alpha2\right) \Gamma 
\left(\frac{\alpha+1}2\right)}{\Gamma(\alpha+1)} 2^\alpha 
[1+(1-\gamma^2)^{1/2}]^{-\alpha}\nonumber\\
&= \frac{\gamma^\alpha}{2\sqrt\pi} \frac{\sqrt{2\pi}\ 2^{-(\alpha-\frac12)} 
\Gamma(\alpha)}{\Gamma(\alpha+1)} 2^\alpha [1+(1-\gamma^2)^{1/2}]^{-\alpha} 
\nonumber\\
\intertext{(by \eqref{eq2.11})}
\label{eq2.12}
&\quad = \alpha^{-1}\gamma^\alpha [1+(1-\gamma^2)^{1/2}]^{-\alpha}.
\end{align}
Then
\[
\int^\infty_0 w^{-1}e^{-w} I_\alpha(\gamma w)dw \le \alpha^{-1} \gamma^\alpha.
\]
Thus part a) is proved.

For part b), we need to differentiate \eqref{eq2.10} under the integral. By 
looking at difference quotients and using the Mean Value Theorem, 
differentiation with respect to $\gamma$ under the integral in the right
hand side of \eqref{eq2.10} 
will be allowed if we can show that for $[a,b] \subseteq (0,1)$,
\begin{equation}\label{eq2.13}
\int^\infty_0 \sup_{\gamma\in [a,b]} |I'_\alpha(\gamma w)| e^{-w} \ dw < 
\infty,
\end{equation}
By formula 8.486.2 on page 970 of Gradshteyn and Ryzhik (1980),
\[
I'_\alpha(z) = \frac12 [I_{\alpha-1}(z) + I_{\alpha+1}(z)].
\]
Hence by \eqref{eq2.7},
\[
|I'_\alpha(z)| \le C(\alpha) e^z[z^{\alpha-1}+z^{\alpha+1}].
\]
In particular,
\[
\sup_{\gamma\in [a,b]} |I'_\alpha(\gamma w)| \le C(\alpha) e^{bw}[w^{\alpha-1} + 
w^{\alpha+1}].
\]
Then since $b<1$ and $\alpha>0$, \eqref{eq2.13} follows.

Thus
\begin{eqnarray*}
\frac{d}{d\gamma} \int^\infty_0 w^{-1} e^{-w}I_\alpha(\gamma w)dw &= &
\int^\infty_0 e^{-w} I'_\alpha(\gamma w)ds\\
&=&\frac12 \int^\infty_0 e^{-w} [I_{\alpha-1} (\gamma w) + I_{\alpha+1}(\gamma 
w)]dw\\
&= &\frac12 \left[\frac{\gamma^{\alpha-1}}{\sqrt{1-\gamma^2}\ 
[1+\sqrt{1-\gamma^2}]^{\alpha-1}}\right]\\
&+&\frac12\left[\frac{\gamma^{\alpha+1}}{\sqrt{1-\gamma^2}\ 
[1 + \sqrt{1-\gamma^2}]^{\alpha+1}}\right]\\
&\le &\frac12 (1-\gamma^2)^{-1/2} [\gamma^{\alpha-1} + \gamma^{\alpha+1}],
\end{eqnarray*}
where we have used  formula 6.611.4 on page 708 of Gradshteyn and Ryzhik (1980) 
for the third equality.  We also see  that the derivative is nonnegative.

To finish, observe that 
\begin{align*}
\left|\frac\partial{\partial\rho} \int^\infty_0 t^{-1} 
e^{-\frac{\rho^2+r^2}{2t}} I_\alpha\left(\frac{\rho r}{t}\right)dt\right| &= 
\left| \frac{d}{d\gamma} \int^\infty_0 w^{-1} e^{-w} I_\alpha(\gamma w)dw\right| 
\left|2r \frac{r^2-\rho^2}{(\rho^2+r^2)^2}\right|\\
&\le r(1-\gamma^2)^{-1/2} (\gamma^{\alpha-1} + \gamma^{\alpha+1}) 
\frac1{\rho^2+r^2},
\end{align*}
as claimed.
\end{proof}

To prove Theorem \ref{thm1.4}, we will need the following consequence of 
\eqref{eq2.10} and \eqref{eq2.12}.

\begin{cor}\label{cor2.4}
For $\gamma = \frac{2\rho r}{\rho^2+r^2} <1$,
$$
\int^\infty_0 t^{-1} e^{-\frac{\rho^2+r^2}{2t}} I_\alpha \left(\frac{\rho 
r}{t}\right) dt = \alpha^{-1} \gamma^\alpha [1+(1-\gamma^2)^{1/2}]^{-\alpha}. 
\eqno \square
$$
\end{cor}

\begin{proof}[Proof of Theorem \ref{thm1.4}]
By Theorem \ref{thm1.3},  \eqref{eq1.7} and Lemma \ref{lem2.1}, for 
$x=\rho\theta$,
\begin{align*}
\frac{d}{dr} P_x(|B_\tau|\le r) &= \int^\infty_0 \intl_{\partial D} \frac12 
\frac\partial{\partial n_y} p_C(t,x,r\eta) r^{n-2} \sin \varphi(\eta) \mu(d\eta) 
dt\\
&= \frac12 r^{\frac{n}2-2} \rho^{1-\frac{n}2} \int^\infty_0 \intl_{\partial D} 
t^{-1} e^{-\frac{\rho^2+r^2}{2t}} \sum^\infty_{j=1}I_{\alpha_j} 
\left(\frac{\rho r}{t}\right) m_j(\theta)\cdot\\
&\quad \cdot\sin\varphi(\eta) 
\frac\partial{\partial n_\eta} m_j(\eta) \mu(d\eta)dt.
\end{align*}

Now for $\gamma = \frac{2\rho r}{\rho^2+r^2} \le 1-\vp$ we have by Corollaries 
\ref{cor2.2} and \ref{cor2.4},
\begin{align*}
&\sum^\infty_{j=1} \int^\infty_0 \intl_{\partial D} t^{-1} 
e^{-\frac{\rho^2+r^2}{2t}} I_{\alpha_j} \left(\frac{\rho r}{t}\right) 
\left|m_j(\theta) \frac\partial{\partial n_\eta} m_j(\eta)\right| \mu(d\eta)dt\\
&\quad \le K \sum^\infty_{j=1} \alpha^{-1}_j \gamma^{\alpha_j} \alpha^{b(n)}_j < 
\infty,
\end{align*}
by \eqref{eq2.9} and the integral test.

Hence we can exchange summation and integration above to get, uniformly 
for $\gamma\le 1-\vp$,
\begin{align*}
\frac{d}{dr} P_x(|B_\tau| \le r) &= \frac12 r^{\frac{n}2-2} \rho^{1-\frac{n}2} 
\sum^\infty_{j=1} \int^\infty_0 t^{-1} e^{-\frac{\rho^2+r^2}{2t}} I_{\alpha_j} 
\left(\frac{\rho r}{t}\right) dt\cdot\\
&\quad\cdot \left[~\intl_{\partial D} \sin \varphi(\eta) \frac\partial{\partial 
n_\eta} m_j(\eta) \mu(d\eta)\right] m_j(\theta)\\
&= \frac12 r^{\frac{n}2-2} \rho^{1-\frac{n}2} \sum^\infty_{j=1} \alpha^{-1}_j 
\gamma^{\alpha_j} [1+(1-\gamma^2)]^{-\alpha_j} \cdot \\ 
&\quad\cdot\left[~\intl_{\partial D} \sin 
\varphi(\eta) \frac\partial{\partial n_\eta} m_j(\eta) \mu(d\eta)\right] 
m_j(\theta),
\end{align*}
as claimed.
\end{proof}

\begin{proof}[Proof of Corollary \ref{cor1.5}]
If $r$ is large,  then $\gamma = \frac{2\rho r}{\rho^2+r^2}$ is small and so by 
Theorem \ref{thm1.4}, as $r\to\infty$
\[
\frac{d}{dr} P_x(|B_\tau|\le r) \sim \frac12 r^{\frac{n}2-2} \rho^{1-\frac{n}2} 
\alpha^{-1}_1 \gamma^{\alpha_1} \left[~\intl_{\partial D} \sin \varphi(\eta) 
\frac\partial{\partial n_\eta} m_1(\eta) \mu(d\eta)\right] m_1(\theta),
\]
where we have used the fact that since
\[
\Delta_{S^{n-1}} m_1 = -\lambda m_1\quad \text{and}\quad m_1\ge 0,
\]
the Hopf Maximum Principle (Protter and Weinberger (1984), Theorem 7 on page 
65) implies $\frac\partial{\partial n_\eta} m_1(\eta)>0$ on $\partial D$. Since 
$\gamma^{\alpha_1} \sim(2\rho)^{\alpha_1} r^{-\alpha_1}$ as $r\to\infty$, we get 
the desired asymptotic upon integrating and appealing to \eqref{eq1.3}.
\end{proof}

\section{Proof of Theorem \ref{thm1.1}}\label{sec3}

\indent

Let $\tau^\pm$ be the first exit times of $X^\pm$ from $C$ and for $u,v>0$, 
define
\[
\eta(-u,v) = \inf\{t>0\colon \ Y_t\notin (-u,v)\}.
\]
For typographical simplicity we write $\tau$ for $\tau_C$. Then for any 
$A\subseteq \partial C$,
\begin{align*}
P_z(Z_\tau \in A)&= P_z(Z_\tau\in A, \text{ exit occurs along $X^-$ path})\\
&\quad + P_z(Z_\tau\in A, \text{ exit occurs along $X^+$ path)}\\
&= P_z(Z_\tau\in A, Y(\eta(-\tau^-,\tau^+)) = -\tau^-)\\
&\quad + P_z(Z_{\tau}\in A, Y(\eta(-\tau^-,\tau^+)) = \tau^+)\\
&= P_z(X^-(\tau^-) \in A, Y(\eta(-\tau^-,\tau^+)) = -\tau^-)\\
&\quad + P_z(X^+(\tau^+) \in A, Y(\eta(-\tau^-,\tau^+)) = \tau^+)\\
&= 2P_z(X^-(\tau^-)\in A, Y(\eta(-\tau^-,\tau^+)) = -\tau^-)
\end{align*}
by independence and symmetry. Writing
\[
f_z(v) = \frac{d}{dv} P_z(\tau^-\le v)
\]
for the density of $\tau^-$, by independence of $X^+$ and $X^-$,
\begin{align}
P_z(Z_\tau\in A) &= 2 \int^\infty_0 P_z(X^-(\tau^-)\in A, Y(\eta(-\tau^-,v) = 
-\tau^-) f_z(v)dv\nonumber\\
&= 2 \int^\infty_0 \intl_{A\times (0,\infty)} P(Y(\eta(-u,v)) = -u) 
P_z((X^-(\tau^-), \tau^-) \in dy \times du) f_z(v)dv\nonumber\\
\label{eq3.1}
&= 2\int^\infty_0 \intl_{A\times (0,\infty)} \frac{v}{u+v} P_z((X^-(\tau^-), 
\tau^-) \in dy\times du)f_z(v)dv.
\end{align}
Hence by Theorem \ref{thm1.3}, Lemma \ref{lem2.1} and \eqref{eq1.7}, for 
$y=r\eta$, $z = \rho\theta$,
\begin{align*}
\frac{d}{dr} P_z(|Z_\tau|\le r) &= \int^\infty_0 \int^\infty_0 \intl_{\partial 
D} \frac{v}{u+v} \frac\partial{\partial n_y} p_C(u,z,y) 
r^{n-2} \sin \varphi(\eta) \mu(d\eta)du \ f_z(v)dv\\
&= r^{\frac{n}2-2} \rho^{1-\frac{n}2} \int^\infty_0 \int^\infty_0 
\intl_{\partial D} \frac{v}{u(u+v)} e^{-\frac{\rho^2+r^2}{2u}} \sum^\infty_{j=1} 
I_{\alpha_j} \left(\frac{\rho r}u\right) m_j(\theta) 
\cdot\\
&\quad \cdot\left[\frac\partial{\partial n_\eta} m_j(\eta)\right]
 \sin \varphi(\eta) \mu(d\eta) du\ f_z(v)dv.
\end{align*}
There is no danger of circular reasoning in using Theorem \ref{thm1.3}  since its 
proof is self-contained. Using Corollary \ref{cor2.2}, if we can show
for $\rho\ne r$,
\begin{equation}\label{eq3.2}
\sum^\infty_{j=1} \alpha^{b(n)}_j \int^\infty_0 \int^\infty_0 \intl_{\partial 
D} \frac{v}{u(u+v)} e^{-\frac{\rho^2+r^2}{2u}} I_{\alpha_j} \left(\frac{\rho 
r}u\right) \mu(d\eta) du\ f_z(v)dv < \infty,
\end{equation}
then by monotone convergence and dominated convergence, exchange of summation 
with integration is allowed and for $\rho \ne r$,
\begin{align}
\frac{d}{dr} P_z(|Z_\tau|\le r) &= r^{\frac{n}2-2} \rho^{1-\frac{n}2} 
\sum^\infty_{j=1} m_j(\theta) \left[~\intl_{\partial D} \sin \varphi(\eta) 
\frac\partial{\partial n_\eta} m_j(\eta) \mu(d\eta)\right]\cdot\nonumber\\
\label{eq3.3}
&\quad  \cdot \int^\infty_0 \int^\infty_0 \frac{v}{u(u+v)} 
e^{-\frac{\rho^2+r^2}{2u}} I_{\alpha_j} \left(\frac{\rho r}u\right) f_z(v)dudv.
\end{align}
The work to justify \eqref{eq3.2} has been done in Section 2:\ The $j^{\rm th}$ 
term is bounded above by
\[
\alpha^{b(n)}_j \mu(\partial D) \int^\infty_0 \frac1u 
e^{-\frac{\rho^2+r^2}{2u}} I_{\alpha_j} \left(\frac{\rho r}{u}\right) du.
\]
Then \eqref{eq3.2} follows from Lemma \ref{lem2.3}a, since $r\ne\rho$.

As it stands, the behavior of $\frac{d}{dr} P_z(|Z_\tau|\le r)$ for large $r$ 
is not obvious from \eqref{eq3.3}. It will turn out that the $j=1$ term dominates. In 
what follows, we write
\begin{equation}\label{eq3.4}
p_j = \alpha_j - \left(\frac{n}2-1\right),
\end{equation}
where $\alpha_j$ is as in Theorem \ref{thm1.4}. 
From \eqref{eq1.2} we have
\begin{equation}\label{eq3.5}
P_z(\tau^->v) \sim C(z) v^{-p_1/2} \quad \text{as}\quad v\to \infty.
\end{equation}

The following lemma will be used to derive asymptotics of the first term in 
\eqref{eq3.3} as well as upper bounds on the remaining terms.

\begin{lem}\label{lem3.1}
Let $\alpha\ge \alpha_1$ and set
\[
I = \int^\infty_0\int^\infty_0 \frac{v}{u(u+v)} e^{-\frac{\rho^2+r^2}{2u}} 
I_\alpha \left(\frac{\rho r}{u}\right) f_z(v) dvdu.
\]
a)~~For some positive $M$ and $K$, independent of $\alpha$,
\[
I\le \begin{cases}
K\alpha^{1/2}(2\rho)^\alpha r^{-2-\alpha},&\text{if $\dfrac{p_1}2 > 1$}\\
\noalign{\smallskip}
K\alpha(2\rho)^\alpha r^{-2-\alpha}\ln r,&\text{if $\dfrac{p_1}2 = 1$}\\
\noalign{\smallskip}
K\alpha^{p_1/2}(2\rho)^\alpha r^{-\alpha-p_1},&\text{if $\dfrac{p_1}2 < 1$,} 
\end{cases}
\]
for $r\ge M$ and $\alpha>1$.
\medskip

\n b)~~For $\frac{p_1}2 >1$,
\[
\lim_{r\to\infty} r^{2+\alpha} I = 2\rho^\alpha E_x(\tau^-).
\]
For $\frac{p_1}2 = 1$,
\[
\lim_{r\to\infty} r^{2+\alpha} (\ln r)^{-1} I = 4\rho^\alpha C(z)
\]
and for $\frac{p_1}2 < 1$,
\[
\lim_{r\to\infty} r^{\alpha+p_1} I = C(z) 2^{p_1/2} \rho^\alpha 
\frac{\Gamma\left(\alpha+\frac{p_1}2\right)}{\Gamma(\alpha+1)} \int^\infty_0 
w^{-p_1/2}(1+w)^{-2}\ dw.
\]
\end{lem}

\begin{proof}
Since
\[
f_z(v) = - \frac{d}{dv} P_z(\tau^->v),
\]
after an integration by parts,
\begin{align}
I &= \int^\infty_0 \int^\infty_0 \frac1{(u+v)^2} P_z(\tau^->v) 
e^{-\frac{\rho^2+r^2}{2u}} I_\alpha \left(\frac{\rho r}{u}\right) \ 
dvdu\nonumber\\
\intertext{(now change of variables $s = \frac{\rho^2+r^2}{2u}$)}
&= \int^\infty_0 \int^\infty_0 \left[\frac{\rho^2+r^2}{2s} + v\right]^{-2} 
P_z(\tau^->v) e^{-s} I_\alpha \left(\frac{2\rho rs}{\rho^2+r^2}\right) 
\frac{\rho^2+r^2}{2s^2}\ dsdv\nonumber\\
&= \frac2{\rho^2+r^2} \int^\infty_0 \int^\infty_0 \left[1 + 
\frac{2sv}{\rho^2+r^2} \right]^{-2} P_z(\tau^->v) e^{-s} I_\alpha 
\left(\frac{2\rho rs}{\rho^2+r^2}\right) \ dsdv\nonumber\\
\label{eq3.5b}
&= \int^\infty_0 \int^\infty_0 H\ dsdv,\quad \text{say.}
\end{align}

\n {\bf Case 1:}\ $\frac{p_1}2 > 1$. By \eqref{eq3.5}, $E_z(\tau^-) < \infty$. 
By \eqref{eq2.7},
\[
H \le K \frac1{\rho^2+r^2}\cdot 1 \cdot P_z(\tau^->v) e^{-s} \left(\frac{\rho 
rs}{\rho^2+r^2}\right)^\alpha \frac{e^\alpha e^{\frac{2\rho rs}{\rho^2+r^2}}}{ 
\left(\alpha + \frac12\right)^\alpha}.
\]
For fixed $\rho>0$, choose $M_1$ independent of $\alpha$ so large that
\begin{equation}\label{eq3.6}
\frac{2\rho r}{\rho^2+r^2} \le \frac12, \qquad r\ge M_1.
\end{equation}
Then 
\begin{equation}\label{eq3.7}
H \le Kr^{-2} P_z(\tau^->v) e^{-s/2} \left(\frac{\rho 
rs}{\rho^2+r^2}\right)^\alpha e^\alpha \left(\alpha + \frac12\right)^{-\alpha}, 
\qquad r\ge M_1.
\end{equation}
Here and in what follows, $K$ will be a number whose exact value might change 
from line to line, but will always be independent of $\alpha$ and $r$. 

First observe  that by \eqref{eq3.7}
\begin{align*}
\int^\infty_0 \int^\infty_0 H\ dsdv &\le \frac{Kr^{-2}e^\alpha}{\left(\alpha + 
\frac12\right)^\alpha} \left(\frac{\rho r}{\rho^2+r^2}\right)^\alpha 
\int^\infty_0 \int^\infty_0 P_z(\tau^->v) e^{-s/2} s^\alpha \ dsdv\\
&= \frac{Kr^{-2}e^\alpha}{\left(\alpha+\frac12\right)^\alpha} \left(\frac{\rho 
r}{\rho^2+r^2}\right)^\alpha E_z(\tau^-) 2^\alpha\Gamma(\alpha+1).
\end{align*}
Now by Stirling's formula
\begin{align*}
\frac{\Gamma(\alpha+1)e^\alpha}{\left(\alpha+\frac12\right)^\alpha} &= 
\frac{e^\alpha\alpha\Gamma(\alpha)}{\left(\alpha + \frac12\right)^\alpha}\\
&\le K \frac{e^\alpha \alpha \alpha^{\alpha-\frac12} e^{-\alpha}}{\left( 
\alpha+\frac12\right)^\alpha}\\
&\le K \alpha^{1/2}.
\end{align*}
Hence for $r\ge M_1$,
\begin{align}
\label{eq3.8}
\int^\infty_0 \int^\infty_0 H\ ds dv &\le Kr^{-2} \left(\frac{2\rho 
r}{\rho^2+r^2}\right)^\alpha \alpha^{1/2}\\
&\le Kr^{-2} \left(\frac{2\rho}r\right)^\alpha \alpha^{1/2}.\nonumber
\end{align}
Referring back to \eqref{eq3.5b}, we see that this gives the desired bound in part a) 
of the lemma.

As for the asymptotic in part b), notice  that 
\[
\int^\infty_0 \int^\infty_0 P_z(\tau^->v) s^\alpha e^{-s/2} \ dsdv < \infty
\]
and
\begin{align*}
\lim_{r\to\infty} r^{2+\alpha} H &= \lim_{r\to\infty} \frac{2r^2}{\rho^2+r^2} 
\left[1 + \frac{2sv}{\rho^2+r^2}\right]^{-2} P_z(\tau^->v) e^{-s} r^\alpha 
I_\alpha 
\left(\frac{2\rho rs}{\rho^2+r^2}\right)\\
&= 2P_z(\tau^->v) e^{-s} \frac{(\rho s)^\alpha}{\Gamma(\alpha+1)},
\end{align*}
using the asymptotic 
\begin{equation}\label{eq3.9}
I_\nu(z) \sim \left(\frac{z}2\right)^\nu/\Gamma(\nu+1)\quad \text{as}\quad z\to 
0.
\end{equation}
Hence by \eqref{eq3.7} and the dominated convergence theorem in \eqref{eq3.5b},
\begin{align*}
\lim_{r\to\infty} r^{2+\alpha} I &= \int^\infty_0 \int^\infty_0 
\lim_{r\to\infty} r^{2+\alpha} H\ ds dv\\
&= \int^\infty_0\int^\infty_0 2P_z(\tau^->v) e^{-s} \frac{(\rho 
s)^\alpha}{\Gamma(\alpha+1)} \ dsdv\\
&= 2\rho^\alpha E_z(\tau^-),
\end{align*}
as claimed.
\medskip

\n {\bf Case 2:}\ $\frac{p_1}2 \le 1$. This part is more delicate because this 
time $E_z(\tau^-) = \infty$. Let $\vp \in\left(0, \frac12\right)$ and use the 
asymptotics \eqref{eq3.5} and \eqref{eq3.9} to choose $M_1$ and $M_2$ such that
\begin{equation}\label{eq3.10}
(1-\vp) C(z) v^{-p_1/2} \le P_z(\tau^->v) \le (1+\vp) C(z) v^{-p_1/2},\qquad 
v\ge M_1
\end{equation}
and
\begin{equation}\label{eq3.11}
\frac{(1-\vp)\left(\frac{z}2\right)^\alpha}{\Gamma(\alpha+1)} \le I_\alpha(z) 
\le \frac{(1+\vp) \left(\frac{z}2\right)^\alpha}{\Gamma(\alpha+1)}, \qquad z\le 
2\rho M_2.
\end{equation}
Notice $M_1$ is independent of $\alpha$ and $M_2$ is not.

We break up the integral $I$ in \eqref{eq3.5b} into three pieces:
\begin{align}
I &= \left(\int^{M_1}_0 \int^\infty_0 + \int^\infty_{M_1} \int^\infty_{M_2r} + 
\int^\infty_{M_1} \int^{M_2r}_0\right) H \ dsdv\nonumber\\
\label{eq3.12}
&= J_1 + J_2 + J_3, \quad \text{say.}
\end{align}
It turns out $J_3$ will dominate as $r\to \infty$.  We have

\begin{align}
J_1 &\le \frac2{r^2} \int^{M_1}_0 \int^\infty_0 1^{-2}\cdot 1\cdot e^{-s} 
I_\alpha \left(\frac{2\rho rs}{\rho^2+r^2}\right) \ dsdv\nonumber\\
&= \frac{2M_1}{r^2}\int^\infty_0 e^{-s} I_\alpha \left(\frac{2\rho 
r}{\rho^2+r^2}s\right)ds\nonumber\\
&= \frac{2M_1}{r^2} \frac{\left(\frac{2\rho 
r}{\rho^2+r^2}\right)^\alpha}{\sqrt{1 - \left(\frac{2\rho 
r}{\rho^2+r^2}\right)^2}\ \left[1 + \sqrt{1 - \left(\frac{2\rho 
r}{\rho^2+r^2}\right)^2}\right]^\alpha}\nonumber\\
\intertext{(Gradshteyn and Ryzhik (1980), 6.611.4 on page 708)}
\label{eq3.13}
&\le K r^{-2-\alpha} (2\rho)^\alpha \quad \text{for}\quad r\ge M_3 \quad  
\text{large,}
\end{align}
where $M_3$ is independent of $\alpha$.

As for $J_2$, by \eqref{eq2.7} 
\begin{align*}
J_2 &\le \frac2{\rho^2+r^2} \int^\infty_{M_1} \int^\infty_{M_2r} 
\left[\frac{2sv}{\rho^2+r^2}\right]^{-2}e^{-s} K \left(\frac{\rho 
rs}{\rho^2+r^2}\right)^\alpha 
\frac{e^\alpha e^{\frac{2\rho 
rs}{\rho^2+r^2}}}{\left(\alpha+\frac12\right)^\alpha} 
\ dsdv\\
&= \frac{K}{\rho^2+r^2} \left(\frac{\rho r}{\rho^2+r^2}\right)^\alpha 
\left[\frac2{\rho^2+r^2}\right]^{-2} 
\frac{e^\alpha}{\left(\alpha+\frac12\right)^\alpha} \int^\infty_{M_2r} 
s^{\alpha-2} e^{-s} e^{\frac{2\rho rs}{\rho^2+r^2}} \ ds\\
&\le K \left(\frac{\rho r}{\rho^2+r^2}\right)^\alpha (\rho^2+r^2) 
\frac{e^\alpha}{\left(\alpha+\frac12\right)^\alpha} \int^\infty_{M_2r} 
s^{\alpha-2} e^{-3s/4} e^{-M_2r/4} e^{\frac{2\rho rs}{\rho^2+r^2}} \ ds\\
&\le K\rho^\alpha r^{2-\alpha} 
\frac{e^\alpha}{\left(\alpha+\frac12\right)^\alpha} e^{-M_2r/4} 
\int^\infty_{M_2r} s^{\alpha-2} e^{-s/2} \ ds,
\end{align*}
for $r$ large, say $r\ge M_4$, where $M_4$ is independent of $\alpha$.
Thus
\begin{align}
J_2&\le K\rho^\alpha r^{2-\alpha} 
\frac{e^\alpha}{\left(\alpha+\frac12\right)^\alpha} e^{-M_2r/4} 2^\alpha 
\Gamma(\alpha-1)\nonumber\\
&\le K(2\rho)^\alpha r^{2-\alpha} e^{-M_2r/4} \frac{e^\alpha}{\left(\alpha + 
\frac12\right)^\alpha} e^{-\alpha+1} (\alpha-1)^{\alpha-\frac32}\nonumber\\
\intertext{(by Stirling's formula)}
\label{eq3.14}
&\le K(2\rho)^\alpha r^{2-\alpha} e^{-M_2r/4} \alpha^{-3/2}.
\end{align}

Now we examine the dominant piece $J_3$. Reversing the order of integration, 
then changing $v$ into $w = \frac{2sv}{\rho^2+r^2}$,
\begin{align}
J_3 &= \int^{M_2r}_0 \int^\infty_{M_1} H\ dvds\qquad \text{(see 
\eqref{eq3.5b})}\nonumber\\
&= \frac2{\rho^2+r^2} \int^{M_2r}_0 \int^\infty_{\frac{2sM_1}{\rho^2+r^2}} 
(1+w)^{-2} P_z\left(\tau^-> \frac{\rho^2+r^2}{2s}w\right) e^{-s} I_\alpha 
\left(\frac{2\rho rs}{\rho^2+r^2}\right) \frac{\rho^2+r^2}{2s}\ dwds\nonumber\\
\label{eq3.15}
&= \int^{M_2r}_0 \int^\infty_{\frac{2sM_1}{\rho^2+r^2}} (1+w)^{-2} P_z 
\left(\tau^- > \frac{\rho^2+r^2}{2s}w\right) e^{-s} s^{-1} I_\alpha 
\left(\frac{2\rho rs}{\rho^2+r^2}\right)\ dwds.
\end{align}
Write
\begin{equation}\label{eq3.16}
h(u) = \int^\infty_u (1+w)^{-2} w^{-p_1/2} \ dw
\end{equation}
and observe
\begin{equation}\label{eq3.17}
\begin{cases}
h(u) \le h(0) < \infty&\text{for $\dfrac{p_1}2 < 1$}\\
\noalign{\smallskip}
h(u) \sim \ln \dfrac1u \text{ as } u\to 0&\text{for $\dfrac{p_1}2 = 1$.}
\end{cases}
\end{equation}
Now for $w > \frac{2sM_1}{\rho^2+r^2}$ we have $\frac{\rho^2+r^2}{2s}w >M_1$, 
hence by \eqref{eq3.10}, \eqref{eq3.15} becomes
\begin{align}
J_3 &\le K \int^{M_2r}_0 \int^\infty_{\frac{2sM_1}{\rho^2+r^2}} (1+w)^{-2} 
\left(\frac{\rho^2+r^2}{2s}w\right)^{-p_1/2} s^{-1} e^{-s} I_\alpha 
\left(\frac{2\rho rs}{\rho^2+r^2}\right)\ dwds\nonumber\\
&= K \left(\frac{\rho^2+r^2}{2}\right)^{-p_1/2} \int^{M_2r}_0 
h\left(\frac{2sM_1}{\rho^2+r^2}\right) s^{p_1/2-1} e^{-s} I_\alpha 
\left(\frac{2\rho rs}{\rho^2+r^2}\right)\ ds\nonumber\\
&\le K(\rho^2+r^2)^{-p_1/2} \int^{M_2r}_0 h\left(\frac{2sM_1}{\rho^2+r^2} 
\right) s^{p_1/2-1} e^{-s} \left(\frac{\rho rs}{\rho^2+r^2}\right)^\alpha 
\frac{e^\alpha e^{\frac{2\rho rs}{\rho^2+r^2}}}{\left(\alpha + 
\frac12\right)^\alpha} \ ds \quad \text{(by \eqref{eq2.7})}\nonumber\\
&= \frac{K(\rho^2+r^2)^{-\alpha-p_1/2} \rho^\alpha r^\alpha e^\alpha}{\left( 
\alpha + \frac12\right)^\alpha} \int^{M_2r}_0 h\left(\frac{2sM_1}{\rho^2+r^2} 
\right) s^{\alpha+p_1/2-1} e^{-s} e^{\frac{2\rho rs}{\rho^2+r^2}}\ ds\nonumber\\
\label{eq3.18}
&\le Kr^{-p_1} \left[\frac{\rho e}{r\left(\alpha+\frac12\right)}\right]^\alpha 
\int^{M_2r}_0 h\left(\frac{2sM_1}{\rho^2+r^2}\right) s^{\alpha+p_1/2-1} e^{-s/2} 
\ ds
\end{align}
for $r$ large, independent of $\alpha$.

If $\frac{p_1}2 < 1$ then by \eqref{eq3.17}, this yields
\begin{align}
J_3 &\le Kr^{-p_1} \left[\frac{\rho e}{r\left(\alpha+\frac12\right)} 
\right]^\alpha 2^\alpha \Gamma \left(\alpha + \frac{p_1}2\right)\nonumber\\
&= Kr^{-p_1} \left[\frac{2\rho e}{r\left(\alpha+\frac12\right)}\right]^\alpha 
\left(\alpha + \frac{p_1}2 - 1\right) \Gamma\left(\alpha + \frac{p_1}2 - 
1\right)\nonumber\\
&\le Kr^{-p_1} \left[\frac{2\rho e}{r\left(\alpha+\frac12\right)}\right]^\alpha 
e^{-\alpha} \left(\alpha + \frac{p_1}2-1\right)^{\alpha+ \frac{p_1}2 -\frac12} 
\qquad \text{(Stirling's formula)}\nonumber\\
\label{eq3.19}
&\le K \alpha^{\frac{p_1-1}2} (2\rho)^\alpha r^{-\alpha-p_1}.
\end{align}
When $\frac{p_1}2 = 1$, observe that for $s\le M_2r$ we have 
$\frac{2sM_1}{\rho^2+r^2} \le \frac{2M_1M_2}r <1$ for $r$ large, independent of 
$\alpha$. Hence by \eqref{eq3.17} and \eqref{eq3.18},
\begin{align}
J_3 &\le \frac{K}{r^2} \left[\frac{\rho e}{r\left(\alpha + \frac12\right)} 
\right]^\alpha \int^{M_2r}_0 \left[\ln \frac{\rho^2+r^2}{2sM_1}\right] s^\alpha 
e^{-s/2} \ ds\nonumber\\
&\le \frac{K}{r^2} \left[\frac{\rho e}{r\left(\alpha+\frac12\right)} 
\right]^\alpha \int^{M_2r}_0 [K\ln r - \ln s] s^\alpha e^{-s/2}\ ds\qquad 
\text{($r$ large)}\nonumber\\
&\le \frac{K}{r^2} \left[\frac{\rho e}{r\left(\alpha + \frac12\right)} 
\right]^\alpha \left[K(\ln r)2^\alpha\Gamma(\alpha+1) - \int^1_0 (\ln s)s^\alpha 
e^{-s/2} \ ds\right]\nonumber\\
&\le \frac{K}{r^2} \left[\frac{\rho e}{r\left(\alpha + \frac12\right)} 
\right]^\alpha \left[K(\ln r) 2^\alpha\Gamma(\alpha+1) - \int^1_0 (\ln s) \ 
ds\right]\nonumber\\
&\le \frac{K}{r^2} \left[\frac{2\rho e}{r\left(\alpha + \frac12\right)} 
\right]^\alpha \Gamma(\alpha+1) \ln r\nonumber\\
&= \frac{K}{r^2} \left[\frac{2\rho e}{r\left(\alpha + \frac12\right)} 
\right]^\alpha \alpha\Gamma(\alpha) \ln r\nonumber\\
\label{eq3.20}
&\le K\alpha(2\rho)^\alpha r^{-\alpha-2} \ln r,\qquad \text{by Stirling's 
formula.}
\end{align}

Combining \eqref{eq3.13}, \eqref{eq3.14}, \eqref{eq3.19} and \eqref{eq3.20}, we 
get the upper bound of part a) in the lemma.

For part b), first assume $\frac{p_1}2 < 1$. For $s\le M_2r$ we have 
$\frac{2\rho rs}{\rho^2+r^2} \le 2\rho M_2$ and for $w > 
\frac{2sM_1}{\rho^2+r^2}$ we have $\frac{\rho^2+r^2}{2s}w > M_1$. Hence by 
\eqref{eq3.10} and \eqref{eq3.11} applied to the $\tau^-$ and $I_\alpha$ factors 
in \eqref{eq3.15}, we get the that integrand in \eqref{eq3.15} is bounded above by
\[
C_\alpha(1+w)^{-2} \left(\frac{\rho^2+r^2}{2s}w\right)^{-p_1/2} s^{-1}e^{-s} 
\left(\frac{\rho rs}{\rho^2+r^2}\right)^\alpha \le C_\alpha r^{-p_1-\alpha} 
w^{-p_1/2} (1+w)^{-2} s^{\alpha+p_1/2-1} e^{-s}.
\]
Moreover, writing \eqref{eq3.15} as
\[
J_3 = \int^\infty_0 \int^\infty_0 F \ dwds,
\]
we see from the asymptotics \eqref{eq3.5} and \eqref{eq3.9} that 
\begin{align*}
\lim_{r\to\infty} r^{p_1+\alpha} F &= \lim_{r\to\infty} r^{p_1+\alpha} 
(1+w)^{-2} C(z) \left(\frac{\rho^2+r^2}{2s}w\right)^{-p_1/2} e^{-s} s^{-1} 
\frac{\left(\frac{\rho rs}{\rho^2+r^2}\right)^\alpha}{\Gamma(\alpha+1)}\\
&= (1+w)^{-2} w^{-p_1/2} C(z) 2^{p_1/2} s^{p_1/2+\alpha-1} e^{-s} \rho^\alpha 
\frac1{\Gamma(\alpha+1)}.
\end{align*}
Since $w^{-p_1/2}(1+w)^{-2} s^{\alpha+p_1/2-1} e^{-s}$ is integrable on 
$s,w>0$, we can apply the dominated convergence theorem to get
\begin{align*}
\lim_{r\to\infty} r^{p_1+\alpha} J_3 &= \int^\infty_0 \int^\infty_0 
\lim_{r\to\infty} r^{p_1+\alpha} F \ dwds\\
&= \frac{C(z) 2^{p_1/2} \rho^\alpha \Gamma\left(\alpha + 
\frac{p_1}2\right)}{\Gamma(\alpha+1)} \int^\infty_0 w^{-p_1/2} (1+w)^{-2} \ dw.
\end{align*}
Combining this with \eqref{eq3.13} and \eqref{eq3.14} and using that 
$\frac{p_1}2 < 1$, we get
\[
\lim_{r\to\infty} r^{p_1+\alpha}I = \lim_{r\to\infty} r^{p_1+\alpha} J_3,
\]
which is the claimed value in part b).

Finally, assume $\frac{p_1}2=1$. Consider the integral
\[
J_4 = \int^{M_2r}_0 \int_{\frac{2sM_1}{\rho^2+r^2}}^\infty (1+w)^{-2} 
\left(\frac{\rho^2+r^2}{2s}w\right)^{-1} e^{-s} s^{-1} 
\left(\frac{\rho rs}{\rho^2+r^2}\right)^\alpha\ dwds
\]
which is just $J_3$ in \eqref{eq3.15} with the factors involving  $\tau^-$ and 
$I_\alpha$  replaced by 
$\left(\frac{\rho^2+r^2}{2s}w\right)^{-1}$ and $\left(\frac{\rho 
rs}{\rho^2+r^2}\right)^\alpha$, respectively. These are more or less the 
asymptotics from \eqref{eq3.5} and \eqref{eq3.9}. Then
\[
J_4 = \int^{M_2r}_0 h\left(\frac{2sM_1}{\rho^2+r^2}\right)  
2(\rho^2+r^2)^{-\alpha-1} (\rho r)^\alpha s^{\alpha} e^{-s}\ ds.
\]
For $s\le M_2r$, $\frac{2sM_1}{\rho^2+r^2} \le \frac{2M_1M_2}r < 1$ for large 
$r$, hence by \eqref{eq3.17}, for such $r$ the integrand of 
$\frac{r^{2+\alpha}}{\ln r} J_4$ 
is bounded above by
\begin{align*}
&C_\alpha \frac{r^{2+\alpha}}{\ln r} \left[\ln \frac{\rho^2+r^2}{2sM_1}\right] 
r^{-\alpha-2} s^{\alpha} e^{-s}\\
\le~&C_\alpha \frac1{\ln r} [K\ln r - \ln s] s^{\alpha} e^{-s}\\
\le~&C_\alpha \frac1{\ln r} [K\ln r + (-\ln s) \vee 0] s^{\alpha} e^{-s}\\
\le~&C_\alpha [K + (-\ln s)\vee 0] s^{\alpha} e^{-s},
\end{align*}
which is integrable on $s>0$. Moreover, the limit of the integrand of 
$r^{2+\alpha}(\ln r)^{-1} J_4$ is
\begin{align*}
&\lim_{r\to\infty} r^{2+\alpha} (\ln r)^{-1} h 
\left(\frac{2sM_1}{\rho^2+r^2}\right) 2(\rho^2+r^2)^{-\alpha-1} (\rho 
r)^\alpha s^{\alpha} e^{-s}\\
&\quad = \lim_{r\to\infty} (\ln r)^{-1} \left[\ln \frac{\rho^2+r^2} 
{2sM_1}\right] 2 \rho^\alpha s^{\alpha} e^{-s}\\
&\quad = 4\rho^\alpha s^{\alpha}e^{-s}.
\end{align*}
Hence by dominated convergence again, 
\[
\lim_{r\to\infty} r^{2+\alpha} (\ln r)^{-1} J_4 = 4\rho^\alpha\Gamma(\alpha+1).
\]
By \eqref{eq3.10}--\eqref{eq3.11}
\[
\frac{(1-\vp)^2C(z)}{\Gamma(\alpha+1)} J_4 \le J_3 \le \frac{(1+\vp)^2 
C(z)}{\Gamma(\alpha+1)} J_4.
\]
Multiply by $r^{2+\alpha}(\ln r)^{-1}$, let $r\to \infty$ then let $\vp\to 0$ to 
end up with
\[
\lim_{r\to\infty} r^{2+\alpha} (\ln r)^{-1} J_3 = 4\rho^\alpha C(z).
\]
By \eqref{eq3.13} and \eqref{eq3.14}, we get
\[
\lim_{r\to\infty} r^{2+\alpha} (\ln r)^{-1} I = 4\rho^\alpha C(z),
\]
as desired.
\end{proof}

Now we can prove Theorem \ref{thm1.1}. Write the sum in \eqref{eq3.3} as 
$\sum\limits^\infty_{j=1} \beta_j(r)$. If we can show that 
\begin{equation}\label{eq3.21}
\sum^\infty_{j=1} \beta_j(r) \sim \beta_1(r) \quad \text{as}\quad r\to\infty,
\end{equation}
then by Lemma \ref{lem3.1} b, the conclusion of Theorem \ref{thm1.1} will hold.
To this end, write
\[
\sum^\infty_{j=1} \beta_j(r) = \beta_1(r) \left[ 1 + \sum^\infty_{j=2} 
\frac{\beta_j(r)}{\beta_1(r)}\right].
\]
It suffices to show
\[
\sum^\infty_{j=2} \frac{\beta_j(r)}{\beta_1(r)} \to 0
\]
as $r\to\infty$. There is no danger in dividing by $\beta_1(r)$ because as in 
the proof of Corollary \ref{cor1.5}, the factor
\[
m_1(\theta) \intl_{\partial D} \sin \varphi(\eta) \frac\partial{\partial n_\eta} 
m_1(\eta) \mu(d\eta)
\]
in $\beta_1$ is positive. By Lemma \ref{lem3.1} a) and \eqref{eq2.6}, for some 
constants $K$ and $M$ independent of $j$, for $r\ge M$
\[
|\beta_j(r)| \le r^{\frac{n}2-2} \rho^{1-\frac{n}2} \begin{cases}
K\alpha^{b(n)+1/2}_j (2\rho)^{\alpha_j} 
r^{-2-\alpha_j},&\text{$\dfrac{p_1}2>1$},\\
\noalign{\smallskip}
K\alpha^{b(n)+1}_j (2\rho)^{\alpha_j} r^{-2-\alpha_j}\ln r,&\text{$\frac{p_1}2 =  
1$,}\\
\noalign{\smallskip}
K\alpha^{b(n)+p_1/2}_j (2\rho)^{\alpha_j} 
r^{-\alpha_j-p_1},&\text{$\dfrac{p_1}2 
< 1.$}
\end{cases}
\]
Then by Lemma \ref{lem3.1} b,
\[
\frac{|\beta_j(r)|}{\beta_1(r)} \le K\alpha^{b(n)+1}_j(2\rho)^{\alpha_j} 
r^{\alpha_1-\alpha_j}.
\]
By the integral test and \eqref{eq2.9}, for any $0<\varepsilon<1$, the series
\[
\sum^\infty_{j=2} \alpha^{b(n)+1}_j \left(\frac{2\rho}r\right)^{\alpha_j}
\]
converges uniformly on $\frac{2\rho}{1-\varepsilon}<r$. Thus, since $\alpha_j > \alpha_1$ for $j\ge 
2$,
\[
\underset{r\to\infty}{\ovl\lim} \left|\sum^\infty_{j=2} 
\frac{\beta_j(r)}{\beta_1(r)}\right| \le \underset{r\to\infty}{\ovl\lim} 
K(2\rho)^{\alpha_1} \sum^\infty_{j=2} \alpha^{b(n)+1}_j 
\left(\frac{2\rho}r\right)^{\alpha_j-\alpha_1} = 0,
\]
and \eqref{eq3.21} follows as desired. $\hfill\square$

\section{Proof of Theorem \ref{thm1.3}}\label{sec4}

\indent

The unbounded, nonsmooth nature of $C$ leads to technicalities not encountered 
in the bounded $C^3$ case considered by Hsu (1986). 
We now state the following result used to prove Theorem \ref{thm1.3}. Before 
giving its proof, we show how it yields Theorem \ref{thm1.3}. We follow Hsu's 
idea of finding the Laplace transform in $t$ of the density. Here and in what 
follows we will write 
\[
B_\vp(y) = \{z\in {\bb R}^n\colon \ |z-y| < \vp\}.
\]

\begin{thm}\label{thm4.1}
a)~~Let $x\in C$, $y\in \partial C\backslash\{0\}$ with $|x|\ne |y|$. Then for 
$\lambda>0$,
\[
\frac\partial{\partial n_y} \int^\infty_0 e^{-\lambda t} p_C(t,x,y)dt = 
\int^\infty_0 e^{-\lambda t} \frac\partial{\partial n_y} p_C(t,x,y) dt.
\]
b)~~For $x\in C$ and $y\in \partial C\backslash\{0\}$, $\frac\partial{\partial 
n_y} p_C(t,x,y) > 0$.\newline
c)~~Let $x\in C$ and $\lambda>0$. If $f\in C^{2,\alpha}(\ovl C)$ is nonnegative 
with compact support in $\ovl C\backslash[\{0\}\cup \partial B_{|x|}(0)]$ then
\[
E_x[e^{-\lambda\tau} f(B_\tau)] = \frac12 \intl_{\partial C} f(y) 
\frac\partial{\partial n_y} G^\lambda_C(x,y) \sigma(dy),
\]
where $G^\lambda_C$ is Green's function for $\frac12 \Delta-\lambda$ on $C$ with 
Dirichlet boundary conditions.$\hfill \square$
\end{thm}

\n {\bf Remarks}
\begin{enumerate}
\item[(1)] Since we use the series expansion of the heat kernel to prove part a), there 
will be an exchange of summation and integration. This requires the condition 
$|x|\ne |y|$ as well as the strange hypothesis about the support of $f$.

\item[(2)]  As we point out below in \eqref{eq4.4}, the function
\[
u(z) = E_z[e^{-\lambda \tau} f(B_\tau)]
\]
solves the boundary value problem
\end{enumerate}

\begin{align*}
&\left(\frac12 \Delta-\lambda\right) u = 0\quad \text{in}\quad C\\
&u|_{\partial C} = f.
\end{align*}
Since a series expansion is known for $p_C(t,x,y)$, the most natural way to try 
to prove part c) of Theorem \ref{thm4.1} is to show directly that
\[
\frac12 \intl_{\partial C} f(y) \frac\partial{\partial n_y} G^\lambda_C(x,y) 
\sigma(dy)
\]
solves the said boundary value problem. It is easy to show the PDE is satisfied, 
but direct verification of the boundary condition eludes us. Hence we are forced 
to take the approach we present below.

\begin{proof}[Proof of Theorem \ref{thm1.3}]
Let $x\in C$ and consider any nonnegative $f\in C^{2,\alpha}(\ovl C)$ with 
compact support in $\ovl C\backslash (\{0\}\cup \partial B_{|x|}(0))$. Then 
since
\begin{equation}\label{eq4.1}
G^\lambda_C(x,y) = \int^\infty_0 e^{-\lambda t} p_C(t,x,y)dt,
\end{equation}
we have
\begin{align*}
&\int^\infty_0 \intl_{\partial C} e^{-\lambda t} f(y) P_x(B_\tau\in dy, \tau\in 
dt) = E_x[e^{-\lambda\tau}f(B_\tau)]\\
&\quad = \frac12 \intl_{\partial C} f(y) \frac\partial{\partial n_y} 
\left[\int^\infty_0 e^{-\lambda t} p_C(t,x,y) dt\right] \sigma(dy)\quad 
\text{(Theorem \ref{thm4.1} c)}\\
&\quad = \frac12 \intl_{\partial C} f(y) \int^\infty_0 e^{-\lambda t} 
\frac\partial{\partial n_y} p_C(t,x,y) dt\ \sigma(dy)\quad \text{(Theorem 
\ref{thm4.1} a)}\\
&\quad = \int^\infty_0 \intl_{\partial C} e^{-\lambda t} f(y) \frac12 
\frac\partial{\partial n_y} p_C(t,x,y) \sigma(dy)dt,
\end{align*}
by Theorem \ref{thm4.1} b and Fubini's Theorem.

Inverting the Laplace transform, we get for any $0\le a<b$,
\[
E_x[f(B_\tau) I_{[a,b]}(\tau)] = \int^b_a \intl_{\partial C} f(y) \frac12 
\frac\partial{\partial n_y} p_C(t,x,y) \sigma(dy) dt.
\]
Varying $f$ appropriately, this yields
\begin{equation}\label{eq4.2}
P_x(B_\tau \in A, a\le\tau\le b) = \int^b_a \intl_A \frac12 
\frac\partial{\partial n_y} p_C(t,x,y) \sigma(dy)dt
\end{equation}
where $A$ is any open subset of $\partial C\backslash(\{0\}\cup \partial 
B_{|x|}(0))$. Since $\partial C\cap \partial B_{|x|}(0)$ is polar with $\sigma$ 
measure 0, and since $\frac\partial{\partial n_y} p_C(t,x,y)$ is continuous as a 
function of $y\in \partial C\backslash\{0\}$ (by Lemma \ref{lem2.1}), we see 
\eqref{eq4.2} holds for arbitrary Borel $A\subseteq \partial C\backslash\{0\}$. 
This yields Theorem \ref{thm1.3}.
\end{proof}

\begin{proof}[Proof of Theorem \ref{thm4.1} a)]
Let $x\in C$, $y\in \partial C\backslash\{0\}$ with $|x|\ne |y|$. Write 
$x=\rho\theta$, $y = r\eta$ and $\gamma = \frac{2\rho r}{\rho^2+r^2}$. Note 
since $\rho\ne r$, $\gamma<1$. Then by Corollary \ref{cor2.2} and Lemma 
\ref{lem2.3} a) the following interchanges of integration, differentiation and 
summation are justified:
\begin{align*}
&\frac\partial{\partial n_y} \int^\infty_0 e^{-\lambda t} p_C(t,x,y)dt = \frac1r 
\frac\partial{\partial n_\eta} \int^\infty_0 e^{-\lambda t} p_C(t,x,y)dt\\
&\quad = \frac1r \frac\partial{\partial n_\eta} \int^\infty_0 e^{-\lambda t} 
t^{-1}(\rho r)^{1-\frac{n}2} e^{-\frac{\rho^2+r^2}{2t}} \sum^\infty_{j=1} 
I_{\alpha_j} \left(\frac{\rho r}t\right) m_j(\theta) m_j(\eta)dt\\
&\quad =\rho^{1-\frac{n}2} r^{-\frac{n}2} \frac\partial{\partial n_\eta} 
\sum^\infty_{j=1} \int^\infty_0 e^{-\lambda t} t^{-1} e^{-\frac{\rho^2+r^2}{2t}} 
I_{\alpha_j} \left(\frac{\rho r}t\right) dt\ m_j(\theta) m_j(\eta)\\
&\qquad \text{(also using monotone and dominated convergence)}\\
&\quad = \rho^{1-\frac{n}2} r^{-\frac{n}2} \sum^\infty_{j=1} \int^\infty_0 
e^{-\lambda t} t^{-1} 
e^{-\frac{\rho^2+r^2}{2t}} I_{\alpha_j} \left(\frac{\rho r}t\right) dt\ 
m_j(\theta) \frac\partial{\partial n_\eta} m_j(\eta)\\
&\quad = \rho^{1-\frac{n}2} r^{-\frac{n}2} \int^\infty_0 e^{-\lambda t} t^{-1} 
e^{-\frac{\rho^2+r^2}{2t}} \sum^\infty_{j=1} I_{\alpha_j} \left(\frac{\rho 
r}t\right) m_j(\theta) \frac\partial{\partial n_\eta} m_j(\eta) dt\\
&\quad =\int^\infty_0 e^{-\lambda t} \frac\partial{\partial n_y} p_C(t,x,y) dt,
\end{align*}
by Lemma \ref{lem2.1}.
\end{proof}

\begin{proof}[Proof of Theorem \ref{thm4.1} b)]
This is an immediate consequence of the Hopf Maximum Principle for parabolic 
operators (Protter and Weinberger (1984), Theorem 6 on page 174).
\end{proof}

Part c) is the hard part. Fix $x\in C$ and write
\begin{equation}\label{eq4.3}
u(z) = E_z[e^{-\lambda\tau} f(B_\tau)],\qquad z\in\ovl C,
\end{equation}
where $f\in C^{2,\alpha}(\ovl C)$ is nonnegative with compact support in $\ovl 
C\backslash(\{0\}\cup \partial B_{|x|}(0))$. Choose $x_0\in C$ such that
\[
|x_0| < |x|\quad \text{and}\quad x_0\notin \text{supp } f.
\]
From now on, $x,f$ and $x_0$ are fixed. Since $\partial C$ is Lipschitz, 
Proposition 8.1.9 and Theorem 8.1.10 on pages 345--346 in Pinsky (1995) imply
\begin{equation}\label{eq4.4}
\left\{\begin{array}{l}
u\in C^{2,\alpha}(C)\\
\noalign{\smallskip}
\left(\dfrac12 \Delta-\lambda\right)u = 0 \text{ in } C\\
\noalign{\smallskip}
u|_{\partial C} = f,\end{array}\right.
\end{equation}
and in fact
\begin{equation}\label{eq4.5}
u(z) = \intl_{\partial C} f(y) k(x;y) \nu_{x_0}(dy)
\end{equation}
where $k(x;y)$ is the Martin kernel with pole at $y\in \partial C$ 
normalized by $k(x_0;y) = 1$ and for any Borel set $A\subseteq \partial C$,
\[
\nu_{x_0}(A) = E_{x_0}[e^{-\lambda \tau}I(B_\tau\in A)].
\]
Note since $f\in C^{2,\alpha}(\ovl C)$, by the Elliptic Regularity Theorem,
\begin{equation}\label{eq4.6}
u\in C^{2,\alpha}(\ovl C\backslash\{0\}).
\end{equation}
Furthermore, since $\partial C$ is Lipschitz, by Theorem 8.1.4 on page 337 in 
Pinsky (1995), any sequence $y_n\in C$ with $y_n\to y\in \partial C$ is a Martin 
sequence. In particular, if $y_n\to y$ along the unit inward normal to $\partial 
C$ at $y\in \partial C\backslash\{0\}$,
\begin{align*}
k(x;y) &= \lim_{n\to\infty} \frac{G^\lambda_C(x,y_n)}{G^\lambda_C(x_0,y_n)}\\
&= \frac{\frac\partial{\partial n_y} G^\lambda_C(x,y)}{\frac\partial{\partial 
n_y} G^\lambda_C(x_0,y)}.
\end{align*}
Hence by \eqref{eq4.5}
\begin{equation}\label{eq4.7}
u(z) = \intl_{\partial C} f(y) \frac\partial{\partial n_y} G^\lambda_C(z,y) 
\nu(dy)
\end{equation}
where
\[
\nu(dy) = \left[\frac\partial{\partial n_y} G^\lambda_C(x_0,y)\right]^{-1} 
\nu_{x_0}(dy).
\]
This representation will allow us to estimate $u$ and its derivatives.

It is known that for $z\in C$,
\begin{equation}\label{eq4.8}
\left(\frac12\Delta-\lambda\right) G^\lambda_C(z,\cdot) = 0\quad \text{on}\quad 
C\backslash\{z\}
\end{equation}
and $G^\lambda_C(z,\cdot)$ is continuous on $\ovl C\backslash\{z\}$ with 
boundary value 0. Hence by the Elliptic Regularity Theorem,
\begin{equation}\label{eq4.9}
G^\lambda_C(z,\cdot) \in C^{2,\alpha}(\ovl C\backslash\{0,z\}).
\end{equation}

To prove part c) of Theorem \ref{thm4.1}, choose $M$ so large and $\vp>0$ so 
small that
\begin{align*}
&\vp < |x_0| < |x|< M\\
&\ovl{B_\vp(0)}\cap \text{supp } f = \emptyset\\
&\ovl{B_M(0)}^c \cap \text{supp} f = \emptyset.
\end{align*}
Then choose $\delta>0$ so small that $\ovl{B_\delta(x)}\subseteq C$. Set
\[
E = C\cap B_M(0) \cap \ovl{B_\delta(x)}^c \backslash \ovl{B_\vp(0)}
\]
(see figure 1).
\medskip

\begin{center}
\epsfig{file=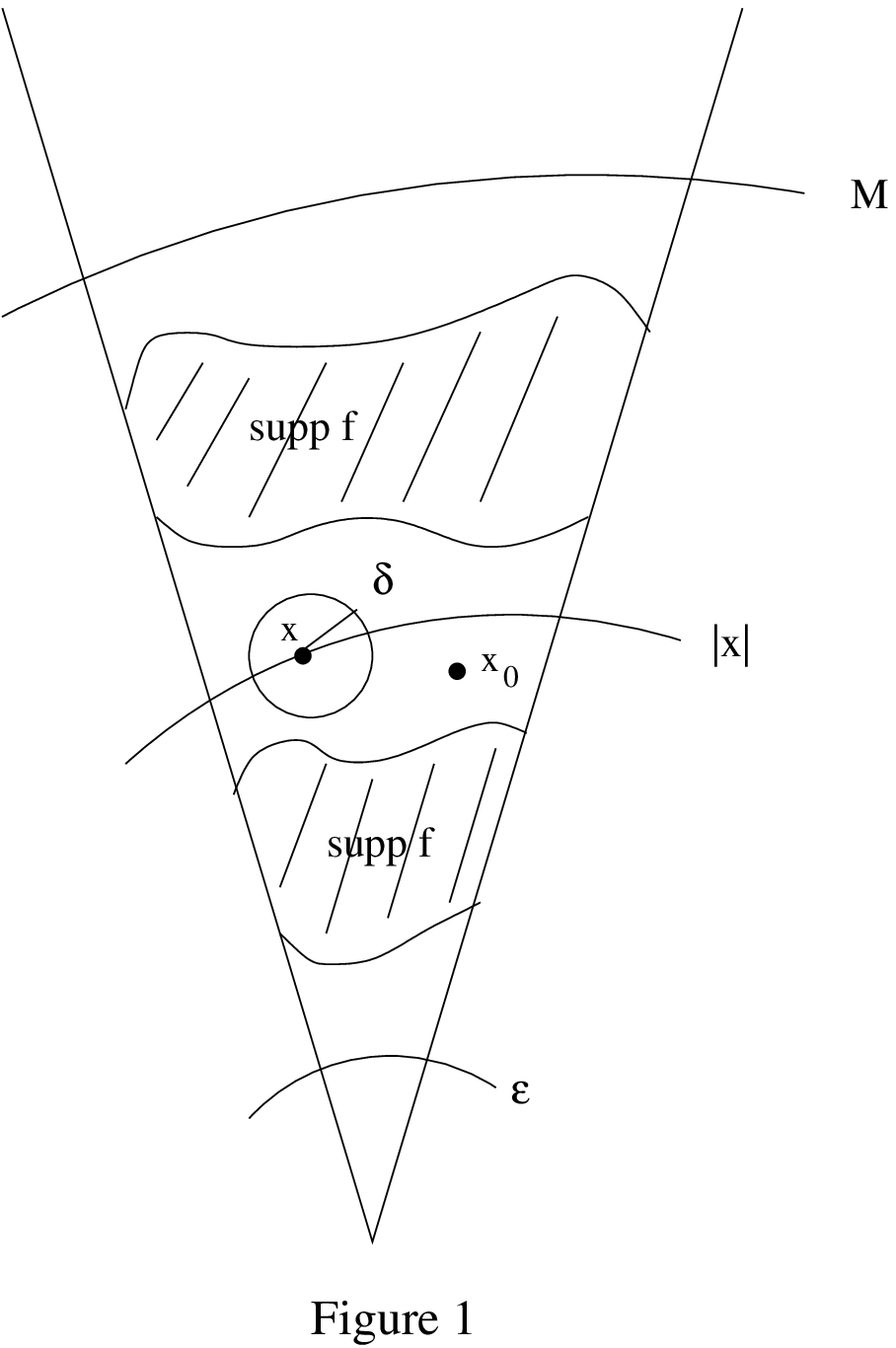,height=4.5in,width=3.3in}
\end{center}

\n By \eqref{eq4.6} and \eqref{eq4.9} we can apply Green's Second Identity:
\begin{align*}
&\intl_E [u(w) \Delta_w G^\lambda_C(x,w) - G^\lambda_C(x,w) \Delta u(w)]dw\\
&\quad = -\intl_{\partial E} \left[u(y) \frac\partial{\partial n_y}  
G^\lambda_C(x,y) - G^\lambda_C(x,y) \frac\partial{\partial n_y} u(y)\right] 
\sigma(dy),
\end{align*}
where $\frac\partial{\partial n_y}$ is the inward normal derivative to $\partial 
E$ and $\sigma(dy)$ is surface  measure on $\partial E$.

By \eqref{eq4.4} and \eqref{eq4.8}, the left-hand side is 0; then breaking up 
$\intl_{\partial E}$ into pieces and solving for the part over $\partial C$,
\begin{align*}
&\intl_{\partial C} \left[u \frac\partial{\partial n_y} G^\lambda_C(x,\cdot) - 
G^\lambda_C(x,\cdot) \frac\partial{\partial n_y} u\right] \sigma(dy)\\
&\quad = \intl_{\partial C\cap B_M(0) \backslash \ovl{B_\delta(0)}}\\
&\quad = -\intl_{\ovl C\cap \partial B_M(0)} + \intl_{\ovl C\cap \partial 
B_\vp(0)} + \intl_{\partial B_\delta(x)}
\end{align*}
(we use the convention that $\frac\partial{\partial n_y}$ is the unit inward 
normal to $\partial B_M(0)$, $\partial B_\vp(0)$, $\partial B_\delta(x)$ 
respectively). Below in Theorem \ref{thm4.7} we will show the first two 
integrals converge to 0 as $M\to \infty$ and $\vp\to 0$. In Theorem \ref{thm4.9} 
we will show the last integral converges to $2u(x)$ as $\delta\to 0$. Thus we 
will end up with
\begin{align*}
u(x) &= \frac12 \intl_{\partial C} \left[u(y) \frac\partial{\partial n_y} 
G^\lambda_C(x,y) - G^\lambda_C(x,y) \frac\partial{\partial n_y} u(y)\right] 
\sigma(dy)\\
&\quad = \intl_{\partial C} f(y) \frac12 \frac\partial{\partial n_y} 
G^\lambda_C(x,y) \sigma(dy),
\end{align*}
using that $u=f$ on $\partial C$ and $G^\lambda_C(x,\cdot) = 0$ on $\partial C$. 
Thus gives part c) of Theorem \ref{thm4.1}.

The representation \eqref{eq4.7} of $u$ will allow us to estimate $u$ and its 
derivatives. For this we need the next result as well as estimates on 
$G^\lambda_C$ and its derivatives.

\begin{lem}\label{lem4.2}
The set $\partial C\cap \text{supp } f$ has finite $\nu$ measure. 
\end{lem}

\begin{proof}
By \eqref{eq4.9}, $G^\lambda_C(x_0,\cdot)\in C^{2,\alpha}(\ovl C\backslash 
\{0,x_0\})$ and by the Hopf Maximum Principle (Protter and Weinberger (1984), 
page 65, Theorem 7),
\[
\frac\partial{\partial n_y} G^\lambda_C(x_0,y)>0,\qquad y\in \partial 
C\backslash 
\{0\}.
\]
Since $\partial C\cap \text{supp } f$ is compact in $\partial C\backslash\{0\}$, 
the desired conclusion follows.
\end{proof}

\begin{lem}\label{lem4.3}
Suppose $\rho>0$, $z\in \ovl C\cap \partial B_\rho(0)$ and $y\in \partial 
C\backslash \{0\}$ with $r = |y|\ne \rho$. If $\gamma = \frac{2\rho 
r}{\rho^2+r^2}$ is positive and sufficiently small, then
\begin{align}
\label{eq4.10}
&\left|\frac\partial{\partial n_y} G^\lambda_C(z,y)\right| \le Kr^{-\frac{n}2} 
\rho^{1-\frac{n}2} \gamma^{\alpha_1}\\
\label{eq4.11}
&\left|\frac\partial{\partial n_z} \frac\partial{\partial n_y} 
G^\lambda_C(z,y)\right| \le Kr^{-\frac{n}2} \left[\rho^{-\frac{n}2} 
\gamma^{\alpha_1} + 
r\rho^{1-\frac{n}2} (\gamma^{\alpha_1-1} + \gamma^{\alpha_1+1}) 
(\rho^2+r^2)^{-1}\right]\\
\label{eq4.12}
&\left|\frac\partial{\partial n_z} G^\lambda_C(z,y)\right| \le 
Kr^{1-\frac{n}2} [\rho^{-\frac{n}2} \gamma^{\alpha_1} + r\rho^{1-\frac{n}2} 
(\gamma^{\alpha_1-1} + \gamma^{\alpha_1+1}) (\rho^2+r^2)^{-1}]\\
\label{eq4.13}
&G^\lambda_C(z,y) \le K(\rho r)^{1-\frac{n}2} \gamma^{\alpha_1},
\end{align}
where $K>0$ is independent of $z$ and $y$.
\end{lem}

\n {\bf Remark}
\begin{enumerate}
\item[(1)] The proof of the bound in \eqref{eq4.13} really only requires $z,y\in \ovl 
C\backslash\{0\}$, $|z|\ne |y|$ and $\frac{2|z||y|}{|z|^2+|y|^2}$ small.

\item[(2)] A similar remark holds for \eqref{eq4.12}. In particular, by the symmetry of 
$G^\lambda_C$, we can replace $\frac\partial{\partial n_z} G^\lambda_C(z,y)$ 
there by $\frac\partial{\partial n_z} G^\lambda_C(y,z)$ and require only that 
$z,y\in \ovl C\backslash\{0\}$ with $|z|\ne |y|$ and 
$\frac{2|z||y|}{|z|^2+|y|^2}$ small.
\end{enumerate}

\begin{proof}
Write $z = \rho\theta$ and $y=r\eta$ in polar coordinates. Then $\rho\ne r$ 
and by Theorem \ref{thm4.1} a), for $\gamma = \frac{2\rho r}{\rho^2+r^2}$,
\begin{align*}
\frac\partial{\partial n_y} G^\lambda_C(z,y) &= \int^\infty_0 e^{-\lambda t} 
\frac\partial{\partial n_y} p_C(t,z,y) dt\\
&\le Kr^{-\frac{n}2} \rho^{1-\frac{n}2} \int^\infty_0 t^{-1} 
e^{-\frac{\rho^2+r^2}{2t}} \sum^\infty_{j=1} I_{\alpha_j} \left(\frac{\rho 
r}t\right) \alpha^{b(n)}_j \ dt\\
\intertext{(by Lemma \ref{lem2.1} and Corollary \ref{cor2.2})}
&\le Kr^{-\frac{n}2} \rho^{1-\frac{n}2} \sum^\infty_{j=1} \int^\infty_0 t^{-1} 
e^{-\frac{\rho^2+r^2}{2t}} I_{\alpha_j} \left(\frac{\rho r}t\right) 
\alpha^{b(n)}_j\ dt\\
\intertext{(Fatou's Lemma)}
&\le Kr^{-\frac{n}2} \rho^{1-\frac{n}2} \sum^\infty_{j=1} \alpha^{b(n)-1}_j 
\gamma^{\alpha_j}\\
\intertext{(by Lemma \ref{lem2.3} a))}
&\le Kr^{-\frac{n}2} \rho^{1-\frac{n}2} \gamma^{\alpha_1} \quad \text{($\gamma$ 
small),}
\end{align*}
which is \eqref{eq4.10}.

Now $\frac\partial{\partial n_z} = -\frac\partial{\partial\rho}$, hence
\begin{align}
\frac\partial{\partial n_z} \frac\partial{\partial n_y} G^\lambda_C(z,y) &= 
-\frac\partial{\partial\rho} 
\int^\infty_0 e^{-\lambda t} r^{-\frac{n}2} \rho^{1-\frac{n}2} t^{-1} 
e^{-\frac{\rho^2+r^2}{2t}} \sum^\infty_{j=1} I_{\alpha_j} \left(\frac{\rho 
r}t\right) m_j(\theta) \frac\partial{\partial n_\eta} m_j(\eta)dt \nonumber\\
\intertext{(by Theorem \ref{thm4.1} a) and Lemma \ref{lem2.1})}
\label{eq4.14}
&= -\frac\partial{\partial\rho} r^{-\frac{n}2} \rho^{1-\frac{n}2} 
\sum^\infty_{j=1} \int^\infty_0 t^{-1} e^{-\lambda t} e^{-\frac{\rho^2+r^2}{2t}} 
I_{\alpha_j} \left(\frac{\rho r}t\right) dt \ m_j(\theta) \frac\partial{\partial 
n_\eta} m_j(\eta)
\end{align}
(by Lemma \ref{lem2.3} and Corollary \ref{cor2.2}). If for $\tilde\gamma = 
\frac{2ur}{u^2+r^2}$ we can show there is $K$ independent of $\alpha,u$ and $r$ 
such that
\begin{align}
&\left|\frac\partial{\partial u} u^{1-\frac{n}2} \int^\infty_0 t^{-1} 
e^{-\lambda 
t} e^{-\frac{u^2+r^2}{2t}} I_\alpha \left(\frac{ur}t\right) dt\right|\nonumber\\
\label{eq4.15}
 &\quad \le 
K[u^{-\frac{n}2} \alpha^{-1} \tilde\gamma^\alpha + 
ru^{1-\frac{n}2} (1-\gamma^2)^{-1/2} (\tilde\gamma^{\alpha-1} + 
\tilde\gamma^{\alpha+1}) (u^2+r^2)^{-1}]
\end{align}
for $0<\tilde\gamma\le\gamma$, then by Corollary \ref{cor2.2} we can 
differentiate under the summation in \eqref{eq4.14}. Moreover, taking $u=\rho$ 
and using Corollary \ref{cor2.2}, for small $\gamma$ we also get the estimate
\begin{align*}
\left|\frac\partial{\partial n_z} \frac\partial{\partial n_y} 
G^\lambda_C(z,y)\right| &\le Kr^{-\frac{n}2} 
\left[\rho^{-\frac{n}2} \sum^\infty_{j=1} \alpha^{b(n)-1}_j \gamma^{\alpha_j} + 
r\rho^{1-\frac{n}2} \sum^\infty_{j=1} \alpha^{b(n)}_j (\gamma^{\alpha_j-1} + 
\gamma^{\alpha_j+1}) (\rho^2+r^2)^{-1}\right]\\
&\le Kr^{-\frac{n}2} [\rho^{-\frac{n}2} \gamma^{\alpha_1} + r\rho^{1-\frac{n}2} 
(\gamma^{\alpha_1-1} + \gamma^{\alpha_1+1}) (\rho^2+r^2)^{-1}],
\end{align*}
giving \eqref{eq4.11}. To prove \eqref{eq4.15}, note by Lemma \ref{lem2.3}
\begin{align*}
&\left|\frac\partial{\partial u} u^{1-\frac{n}2} \int^\infty_0 t^{-1} 
e^{-\lambda t} e^{-\frac{u^2+r^2}{2t}} I_\alpha\left(\frac{ur}t\right) 
dt\right|\\
&\quad \le \left(1-\frac{n}2\right) u^{-\frac{n}2} K\alpha^{-1} \tilde 
\gamma^\alpha + u^{1-\frac{n}2} r(1-\tilde\gamma^2)^{-1/2} 
(\tilde\gamma^{\alpha-1} + \tilde\gamma^{\alpha+1}) (u^2+r^2)^{-1}\\
&\quad \le K[u^{-\frac{n}2} \alpha^{-1} \tilde\gamma^\alpha + ru^{1-\frac{n}2} 
(1-\gamma^2)^{-1/2} (\tilde\gamma^{\alpha-1} + \tilde\gamma^{\alpha+1}) 
(u^2+r^2)^{-1}]
\end{align*}
as desired

To  prove \eqref{eq4.12}, we repeat the proof of \eqref{eq4.11} almost word for 
word. The only change is in \eqref{eq4.14} where the initial $r^{-\frac{n}2}$ is 
replaced by $r^{1-\frac{n}2}$ and the $\frac\partial{\partial n_\eta} m_j(\eta)$ 
is replaced by $m_j(\eta)$.

For the proof of \eqref{eq4.13}, note by \eqref{eq4.1}, \eqref{eq2.1} and 
Corollary \ref{cor2.2},
\begin{align*}
G^\lambda_C(z,y) &\le (\rho r)^{1-\frac{n}2} \sum^\infty_{j=1} \int^\infty_0 
t^{-1} 
e^{-\frac{\rho^2+r^2}{2t}} I_{\alpha_j} \left(\frac{\rho r}t\right) dt \ 
\alpha^{b(n)}_j\\
&\le K(\rho r)^{1-\frac{n}2} \sum^\infty_{j=1} \alpha^{b(n)-1}_j \gamma^{\alpha_j} 
\quad \text{(by Lemma \ref{lem2.3} a))}\\
&\le K(\rho r)^{1-\frac{n}2} \gamma^{\alpha_1},\quad \text{for $\gamma$ small.} 
\qquad \qed
\end{align*}
\renewcommand{\qed}{}\end{proof}

\begin{cor}\label{cor4.4}
Given $z\in \ovl C\backslash\{0\}$, for any compact set $E\subseteq \partial 
C\backslash\{0\}$ with $E\cap \partial B_{|z|}(0) = \emptyset$, there is a 
neighborhood $N$ of $z$ in $\ovl C\backslash\{0\}$ such that
\[
\sup\left\{\left|\nabla_w\frac\partial{\partial n_y} 
G^\lambda_C(w,y)\right|\colon \ w\in N, y\in E\right\}< \infty.
\]
\end{cor}

\begin{proof}
For $\gamma = \frac{2|w||y|}{|w|^2+|y|^2}$, we have for some neighborhood $N$ of 
$z$ in $\ovl C\backslash\{0\}$, $\sup \{\gamma\colon \ w\in N,y\in E\} < 1$. 
Then we can use Lemma \ref{lem2.3} and Corollary \ref{cor2.2} as we did in the 
proof of \eqref{eq4.11} above to get the desired conclusion.
\end{proof}

\begin{cor}\label{cor4.5}
Suppose $\rho \ne |x|$ and $z\in \ovl C\cap \partial B_\rho(0)$. Then for some 
constant $K$ independent of $z$,
\begin{align*}
G^\lambda_C(x,z) &\le K \begin{cases}
\rho^{1-\frac{n}2 + \alpha_1},&\text{$\rho$ small,}\\
\rho^{1-\frac{n}2-\alpha_1},&\text{$\rho$ large,}\end{cases}\\
\left|\frac\partial{\partial n_z} G^\lambda_C(x,z)\right| &\le K \begin{cases}
\rho^{-\frac{n}2 +\alpha_1},&\text{$\rho$ small,}\\
\rho^{-\frac{n}2 - \alpha_1},&\text{$\rho$ large,}\end{cases}
\end{align*}
where $\frac\partial{\partial n_z}$ is the inward normal derivative on $\partial 
B_\rho(0)$.
\end{cor}

\begin{proof}
If $\rho$ is small or large, $\gamma = \frac{2\rho|x|}{\rho^2+|x|^2}$ is small. 
Hence by the remark after Lemma \ref{lem4.3},
\[
G^\lambda_C(x,z) \le K|x|^{1-\frac{n}2} \rho^{1-\frac{n}2} \gamma^{\alpha_1}
\]
and
\[
\left|\frac\partial{\partial n_z} G^\lambda_C(x,z)\right| \le K|x|^{1-\frac{n}2} 
[\rho^{-\frac{n}2} \gamma^{\alpha_1} + |x|\rho^{1-\frac{n}2} 
(\gamma^{\alpha_1-1} + \gamma^{\alpha_1+1}) (\rho^2+|x|^2)^{-1}].
\]
The desired bounds follow from these inequalities.
\end{proof}

As another application of Lemma \ref{lem4.3}, we get bounds on the function
\[
u(z) =  \intl_{\partial C} f(y) \frac\partial{\partial n_y} 
G^\lambda_C(z,y) \nu(dy)
\]
(from \eqref{eq4.7}) and its normal derivatives for large and small $z$.

\begin{lem}\label{lem4.6}
a)~~For some constant $K$, for $z\in\ovl C\backslash\{0\}$ with $|z|=M$ 
sufficiently large
\begin{align*}
|u(z)| &\le KM^{1-\frac{n}2-\alpha_1}\\
\left|\frac\partial{\partial n_z}u(z)\right| &\le KM^{-\frac{n}2-\alpha_1}
\end{align*}
where $\frac\partial{\partial n_z}$ is the inward normal derivative on $\partial 
B_M(0)$.

b)~~For some constant $K$, for $z\in\ovl C\backslash\{0\}$ with $|z|=\vp$ 
sufficiently small,
\begin{align*}
|u(z)| &\le K \vp^{1-\frac{n}2+\alpha_1}\\
\left|\frac\partial{\partial n_z} u(z)\right| &\le K\vp^{-\frac{n}2+ \alpha_1}.
\end{align*}
\end{lem}

\begin{proof}
Note if $y\in \partial C\cap \text{supp } f$, then for $\rho$ large or $\rho$ 
small, $\frac{2\rho |y|}{\rho^2+|y|^2}$ is uniformly small in $y$. Then by 
\eqref{eq4.7} and \eqref{eq4.10} with $|z|=\rho = M$ or $\vp$,
\[
|u(z)| \le K \intl_{\partial C} f(y) |y|^{-\frac{n}2} \rho^{1-\frac{n}2} 
\left(\frac{2\rho|y|}{\rho^2+|y|^2}\right)^{\alpha_1} \nu(dy).
\]
Since $\text{supp } f$ is compact in $\ovl C\backslash(\{0\}\cup \partial 
B_{|x|}(0))$, we have $0<\inf\{|y|\colon \ y\in \text{supp } f\} \le 
\sup\{|y|\colon \ y\in \text{supp } f\} < \infty$ and so by Lemma \ref{lem4.2}, 
when $\rho =\vp$,
\[
|u(z)| \le K\vp^{1-\frac{n}2+\alpha_1},\quad \vp \text{ small;}
\]
and when $\rho =M$,
\[
|u(z)| \le KM^{1-\frac{n}2-\alpha_1},\quad M \text{ large.}
\]

The derivative estimates are a bit more delicate. Looking at the last term in 
\eqref{eq4.11}:\ for $y\in \text{supp } f$,
\begin{align*}
\rho^{1-\frac{n}2} (\gamma^{\alpha_1-1} + \gamma^{\alpha_1+1}) (\rho^2+r^2)^{-1} 
&= \rho^{1-\frac{n}2} \left[\left(\frac{2\rho r}{\rho^2+r^2}\right)^{\alpha_1-1} 
 + \left(\frac{2\rho 
r}{\rho^2+r^2}\right)^{\alpha_1+1}\right] (\rho^2+r^2)^{-1}\\
&\le K\rho^{1-\frac{n}2} \begin{cases}
[\rho^{1-\alpha_1} + \rho^{-1-\alpha_1}]\rho^{-2},&\text{$\rho$ large,}\\
\rho^{\alpha_1-1} + \rho^{\alpha_1+1},&\text{$\rho$ small,}\end{cases}\\
&\le K \begin{cases}
\rho^{-\alpha_1-\frac{n}2},&\text{$\rho$ large,}\\
\rho^{\alpha_1-\frac{n}2},&\text{$\rho$ small.}\end{cases}
\end{align*}
Thus for $y\in \text{supp } f$ and $|z| = M$ large, \eqref{eq4.11} yields
\begin{align*}
\left|\frac\partial{\partial n_z} \frac\partial{\partial n_y} 
G^\lambda_C(z,y)\right| &\le K\left[M^{-\frac{n}2} 
\left(\frac{2Mr}{M^2+r^2}\right)^{\alpha_1} + M^{-\alpha_1-\frac{n}2}\right]\\
&\le K [M^{-\frac{n}2-\alpha_1} + M^{-\alpha_1-\frac{n}2}]\\
&\le KM^{-\alpha_1-\frac{n}2}
\end{align*}
and if $|z|=\vp$ is small,
\begin{align*}
\left|\frac\partial{\partial n_z} \frac\partial{\partial n_y} 
G^\lambda_C(z,y)\right| &\le K\left[\vp^{-\frac{n}2} \left(\frac{2\vp 
r}{\vp^2+r^2}\right)^{\alpha_1} + \vp^{\alpha_1-\frac{n}2}\right]\\
&\le K[\vp^{-\frac{n}2 +\alpha_1} + \vp^{\alpha_1-\frac{n}2}]\\
&= K\vp^{-\frac{n}2 +\alpha_1}.
\end{align*}
Then if we can differentiate under the integral 
\begin{align*}
\left|\frac\partial{\partial n_z} u(z)\right| &\le  \intl_{\partial C} 
f(y) \left|\frac\partial{\partial n_z} \frac\partial{\partial n_y} 
G^\lambda_C(z,y)\right| \nu(dy)\\
&\le \begin{cases}
KM^{-\alpha_1-\frac{n}2},&\text{if $|z|=M$ is large}\\
K\vp^{-\frac{n}2+\alpha_1},&\text{if $|z|=\vp$ is small}
\end{cases}
\end{align*}
(by Lemma \ref{lem4.2}), as desired.

As for differentiation under the integral, by bounding difference quotients via 
the Mean Value Theorem, it is easy to see by Corollary \ref{cor4.4} the exchange 
is justified.
\end{proof}

\begin{thm}\label{thm4.7}
We have
\[
\intl_{\ovl C\cap \partial B_M(0)} \left[u(z) \frac\partial{\partial n_z} 
G^\lambda_C(x,z) - G^\lambda_C(x,z) \frac\partial{\partial n_z} u(z)\right] 
\sigma(dz) \to 0
\]
as $M\to \infty$ or as $M\to 0$. Here $\frac\partial{\partial n_z}$ is the unit 
inward normal derivative on $\partial B_M(0)$.
\end{thm}

\begin{proof}
By Corollary \ref{cor4.5} and Lemma \ref{lem4.6},
\begin{align*}
\text{integrand} &\le K \begin{cases}
M^{1-\frac{n}2+\alpha_1} \cdot M^{-\frac{n}2+\alpha_1} + M^{1-\frac{n}2 + 
\alpha_1} \cdot M^{-\frac{n}2+\alpha_1},&\text{$M$ small,}\\
M^{1-\frac{n}2 - \alpha_1}\cdot M^ {-\frac{n}2-\alpha_1} + M^{1-\frac{n}2 - 
\alpha_1} \cdot M^{-\frac{n}2-\alpha_1},&\text{$M$ large,}\end{cases}
\\
&= K \begin{cases}
M^{1-n+2\alpha_1},&\text{$M$ small,}\\
M^{1-n-2\alpha_1},&\text{$M$ large.}\end{cases}
\end{align*}
Since $\sigma(\partial B_M(0)) \le KM^{n-1}$, we get the desired conclusion.
\end{proof}

The next order of business is to study $G^\lambda_C(x,\cdot)$ in a small 
neighborhood of $x$. To this end, introduce the function
\begin{equation}\label{eq4.16}
G^\lambda(x,y) := \int^\infty_0 e^{-\lambda t} p(t,x,y) dt
\end{equation}
where
\begin{equation}\label{eq4.17}
p(t,x,y) = (2\pi t)^{-n/2} \exp \left(-\frac1{2t} |x-y|^2\right)
\end{equation}
is the usual Gaussian kernel. The relevant properties of $G^\lambda$ are stated 
in the next lemma.

\begin{lem}\label{lem4.8}
a)~~For $0<|x-y|$ small,
\[
G^\lambda(x,y) \le K \begin{cases}
|x-y|^{2-n},&\text{$n\ge 3$,}\\
-\ln|x-y|,&\text{$n=2$.}\end{cases}
\]

b)~~As $|x-y|\to 0$,
\[
\frac{y-x}{|y-x|} \cdot \nabla_y G^\lambda(x,y) \sim -\pi^{-n/2} 
\Gamma\left(\frac{n}2\right) |x-y|^{1-n}.
\]

c)~~For some small neighborhood $N$ of $x$ with compact closure in $C$,
\[
\sup\{|\nabla_y G^\lambda(w,y)|\colon \ w\in \partial C\backslash\{0\}, y\in N\} 
< \infty.
\]
\end{lem}

\begin{proof}
After changing variables $u=\lambda t$,
\begin{align*}
G^\lambda(x,y) &= \int^\infty_0 e^{-\lambda t} (2\pi t)^{-n/2} \exp\left(- 
\frac1{2t} |x-y|^2\right)dt\\
&= (2\pi)^{-\frac{n}2} \lambda^{\frac{n}2-1} \int^\infty_0 u^{-\frac{n}2}e^{-u} 
\exp\left(-\frac\lambda2 \frac{|x-y|^2}u\right) du\\
&= 2^{\frac12(1-\frac{n}2)} \pi^{-\frac{n}2} \lambda^{\frac12(\frac{n}2-1)} 
|x-y|^{1-\frac{n}2} K_{\frac{n}2-1} (\sqrt{2\lambda}\ |x-y|)
\end{align*}
(by formula 3.471.12 on page 340 in Gradshteyn and Ryzhik (1980)) where $K_\nu$ is 
the modified Bessel function. It is known that
\begin{align}
\label{eq4.18}
K_0(z) &\sim - \ln z \quad \text{as}\quad z\to 0\\
\label{eq4.19}
K_\nu(z) &\sim \frac12 \Gamma(\nu) \left(\frac{z}2\right)^{-\nu} \quad\text{as} 
\quad z\to 0\\
\label{eq4.20}
K'_\nu(z) &\sim -\sqrt{\frac\pi{2z}}\ e^{-z}\quad \text{as}\quad z\to\infty\\
\label{eq4.21}
K_\nu(z) &\sim \sqrt{\frac\pi{2z}} e^{-z} \quad \text{as}\quad z\to\infty\\
\label{eq4.22}
K'_\nu(z) &= -K_{\nu+1}(z) + \frac\nu{z} K_\nu(z)
\end{align}
(\eqref{eq4.18}--\eqref{eq4.21} can be found in Abramowitz and Stegun (1972) 
pp.\ 375--378, formulas 9.6.8, 9.6.9, 9.7.4, 9.7.2, respectively. Formula 
\eqref{eq4.22} is from Watson (1922) page 79, formula (4) in section 3.71).

 Part a) is an immediate consequence of 
\eqref{eq4.18}--\eqref{eq4.19}. 

Now for $C_1 = 2^{\frac12(1-\frac{n}2)} \pi^{-\frac{n}2} 
\lambda^{\frac12(\frac{n}2-1)}$, by \eqref{eq4.22} and \eqref{eq4.19}, as 
$|x-y|\to 0$,
\begin{align}
\frac{y-x}{|y-x|} \cdot \nabla_y G^\lambda(x,y) &= 
C_1\left[\left(1-\frac{n}2\right) |x-y|^{-\frac{n}2} 
K_{\frac{n}2-1} (\sqrt{2\lambda} |x-y|)\right.\nonumber\\
&\quad\left. + |x-y|^{1-\frac{n}2} \sqrt{2\lambda}\ K'_{\frac{n}2-1} 
(\sqrt{2\lambda} |x-y|)\right]\nonumber\\
&= -C_1|x-y|^{1-\frac{n}2} \sqrt{2\lambda}\ K_{\frac{n}2}(\sqrt{2\lambda}\ 
|x-y|)\nonumber\\
&\sim -C_1 |x-y|^{1-\frac{n}2} \sqrt{2\lambda}\ \frac12 
\Gamma\left(\frac{n}2\right) \left(\sqrt{\frac\lambda2}\ 
|x-y|\right)^{-\frac{n}2}\nonumber\\
&= -\pi^{-n/2} \Gamma\left(\frac{n}2\right) |x-y|^{1-n}.\nonumber
\end{align}
This gives part b).

Finally, for some $\delta>0$, $\ovl{B_\delta(x)} \subseteq C$ and so 
$\inf\{|w-y|\colon \ w\in \partial C\backslash\{0\}$, $y\in B_\delta(x)\} > 0$. 
Then by \eqref{eq4.20}--\eqref{eq4.21}, part c) follows.
\end{proof}

\begin{thm}\label{thm4.9}
We have
\[
\lim_{\delta\to 0} \intl_{\partial B_\delta(x)} \left[u(z) 
\frac\partial{\partial n_z} G^\lambda_C(x,z) - G^\lambda_C(x,z) 
\frac\partial{\partial n_z} u(z)\right] \sigma(dz) = 2u(x)
\]
where $\frac\partial{\partial n_z}$ is the unit inward normal derivative on 
$\partial B_\delta(0)$.
\end{thm}

\begin{proof}
By \eqref{eq4.6} $\frac\partial{\partial n_z}u$ is bounded on a neighborhood of 
$x$. Since $G^\lambda_C\le G^\lambda$, by Lemma \ref{lem4.8} a) for $\delta$ 
small,
\[
\intl_{\partial B_\delta(x)} \left|G^\lambda_C(x,z)  \frac\partial{\partial n_z} 
u(z)\right| \sigma(dz) \le K\sigma(\partial B_\delta(0)) 
\left\{\begin{array}{ll}
\delta^{2-n},&n\ge 3\\
-\ln \delta,&n=2\end{array}\right\} \to 0,
\]
as  $\delta\to 0$.
Thus we just need to show
\begin{equation}\label{eq4.23}
\lim_{\delta\to  0} \intl_{\partial B_\delta(0)} \left[u(z) 
\frac\partial{\partial n_z} G^\lambda_C(x,z)\right] \sigma(dz) = 2u(x).
\end{equation}
It is well-known that for the exit time $\tau$ of Brownian motion $B_t$ from 
$C$,
\[
p_C(t,x,z) = p(t,x,z)  -E_x [I_{\tau<t} p(t-\tau, B_\tau,z)].
\]
Then
\[
G^\lambda_C(x,z) = G^\lambda(x,z) - E_x[e^{-\lambda\tau} G^\lambda(B_\tau,z)],
\]
and as a consequence, for $z\in \partial B_\delta(x)$
\[
\frac\partial{\partial n_z} G^\lambda_C(x,z) = \frac\partial{\partial n_z} 
G^\lambda(x,z) - E_x \left[e^{-\lambda\tau} \frac\partial{\partial n_z} 
G^\lambda(B_\tau,z)\right].
\]
The exchange of $\frac\partial{\partial n_z}$ and $E_x$ is justified as follows. 
Bound difference quotients via the Mean Value Theorem. Then the exchange is 
justified provided
\[
\sup\{|\nabla_z G^\lambda(w,z)|\colon \ z\in B_{2\delta}(x), w\in \partial 
C\backslash\{0\}\} < \infty.
\]
This follows from Lemma \ref{lem4.8}, part c.

Furthermore, since $u$ is bounded near $x$, by Lemma \ref{lem4.8} c
\[
\left|~\intl_{\partial B_\delta(x)} u(z) E_x\left[e^{-\lambda\tau} 
\frac\partial{\partial n_z} G^\lambda(B_\tau,z)\right] \sigma(dz)\right| \le 
K\sigma(\partial B_\delta(x)) \to 0 \text{ as } \delta \to 0.
\]
Thus to get \eqref{eq4.23} we need only show
\[
\lim_{\delta\to 0} \intl_{\partial B_\delta(x)} u(z) \frac\partial{\partial n_z} 
G^\lambda(x,z) \sigma(dz) = 2u(x).
\]
By Lemma \ref{lem4.8} b, on $\partial B_\delta(x)$, $\frac\partial{\partial 
n_z} G^\lambda(x,z) = -\frac{z-x}{|z-x|} \cdot \nabla_zG^\lambda(x,z) \sim 
\pi^{-n/2} \Gamma\left(\frac{n}2\right) \delta^{1-n}$ 
as $\delta\to 0$. Since $u$ is continuous and bounded near $x$,
\begin{align*}
\lim_{\delta\to 0} \intl_{\partial B_\delta(x)} u(z) \frac\partial{\partial n_z} 
G^\lambda(x,z) \sigma(dz) &= u(x) \pi^{-n/2} \Gamma\left(\frac{n}2\right) 
\sigma(\partial B_1(x))\\
&= 2u(x),
\end{align*}
as desired.
\end{proof}

\end{document}